\documentclass{article}
\usepackage[latin1]{inputenc}
\usepackage[T1]{fontenc}
\usepackage{courier}
\usepackage[english]{babel}
\usepackage{listings}
\usepackage{amssymb}
\usepackage{amsmath}
\usepackage{amsthm}
\usepackage{bm}
\usepackage{graphicx}
\usepackage[margin=3cm]{geometry}
\usepackage{url}
\usepackage[font=small,labelfont=bf]{caption}
\usepackage{epstopdf}

\usepackage{enumitem}
\usepackage{color}

\usepackage{algorithm}
\usepackage{algpseudocode}

\newcommand{\im}{x}
\newcommand{\imorig}{\im_\text{orig}}
\newcommand{\imsol}{\im^*}
\newcommand{\sysmat}{A}
\newcommand{\sino}{b}

\newcommand{\ellone}{\text{P$_1$}}
\newcommand{\ellp}{\text{LP}}
\newcommand{\tv}{\text{TV}}
\newcommand{\nv}{N_\text{v}}
\newcommand{\nside}{N_\text{side}}

\newcommand{\imsep}{\qquad}

\newcommand{\sm}[1]{\mathit{#1}}

\DeclareMathOperator*{\argmin}{arg\,min}
\begin{document}
\title{How little data is enough? Phase-diagram analysis of sparsity-regularized X-ray CT}
\author{Jakob S. J\o{}rgensen\footnote{Corresponding author. Email: jakj@dtu.dk. Technical University of Denmark, Department of Applied Mathematics and Computer Science, Richard Petersens Plads, 2800 Kongens Lyngby, Denmark.} ~and Emil Y. Sidky\footnote{University of Chicago, Department of Radiology MC-2026, 5841 South Maryland Avenue, Chicago, IL 60637, USA.}}

 \maketitle
 
 \begin{abstract}
We introduce phase-diagram analysis, a standard tool in compressed sensing, to the X-ray CT community as a systematic method for determining how few projections suffice for accurate sparsity-regularized reconstruction. 
In compressed sensing a phase diagram is a convenient way to study and express certain theoretical relations between sparsity and sufficient sampling. 
We adapt phase-diagram analysis for empirical use in X-ray CT for which the same theoretical results do not hold. We demonstrate in three case studies the potential of phase-diagram analysis for providing quantitative answers to questions of undersampling: First we demonstrate that there are cases where X-ray CT empirically performs comparable with an optimal compressed sensing strategy, namely taking measurements with Gaussian sensing matrices. Second, we show that, in contrast to what might have been anticipated, taking randomized CT measurements does not lead to improved performance compared to standard structured sampling patterns. Finally, we show preliminary results of how well phase-diagram analysis can predict the sufficient number of projections for accurately reconstructing a large-scale image of a given sparsity by means of total-variation regularization.
\end{abstract}
 
 \section{Introduction} \label{sec:introduction}
\subsection{Sparsity regularization in X-ray CT}
Sparsity-regularized (SR) image reconstruction has shown great promise for X-ray CT. Many works, e.g, \cite{SidkyTV:06,song2007sparseness,sidky2008image,chen2008prior,Bian:10,ritschl2011improved}, have demonstrated that accurate reconstructions
can be obtained from substantially less projection data than is normally required by standard analytical methods such as filtered back-projection and algebraic reconstruction methods.
Acquiring less data is of interest in many applications of X-ray CT to reduce scan time or exposure to ionizing radiation.

The typical SR setup for X-ray CT, and the one we employ, is that an unknown discrete image $\im \in \mathbb{R}^N$ is to be reconstructed from measured discrete data $\sino \in \mathbb{R}^m$, connected to $\im$ through a linear model, $\sino \approx \sysmat \im$, for some measurement matrix $\sysmat \in \mathbb{R}^{m\times{}N}$. A common reconstruction problem is 
\begin{equation}
 \imsol = \argmin_\im ~R(\im) \quad \text{subject to} \quad  \|\sysmat \im - \sino\|_2 \leq \epsilon, \label{eq:ineqreconprob}
\end{equation}
where $R(\im)$ is a sparsity regularizer, for example the $1$-norm, the total variation (TV) semi-norm, or a $1$-norm of wavelet coefficients or coefficients in a learned dictionary, depending on which domain sparsity is expected in, and $\epsilon$ is a regularization parameter that must be chosen to balance the level of regularization enforced with the misfit to data.

In contrast to analytical and algebraic reconstruction methods, SR can admit reconstructions in the underdetermined case $m < N$ as shown in the references given above. However, from the existing individual studies it is difficult to synthesize a coherent quantitative understanding of the undersampling potential of SR in CT. 
From a practical point of view, we want to know how many CT projections to acquire in order to obtain a SR reconstruction of sufficient quality to reliably solve the relevant imaging task, for example detection, classification, segmentation, etc. This question is difficult to address meaningfully in general, because specific applications pose different challenges, for example varying levels of noise and inconsistencies in the data as well as different quality requirements on the reconstruction. 
But even in an application-independent setting, systematic analysis of the undersampling potential of SR in CT remains unexplored.

We consider in the present work an idealized form of the reconstruction problem \eqref{eq:ineqreconprob} with $\epsilon = 0$ and consider only synthetic noise-free data. This simplified setup allows us to study more precise questions with fewer complicating factors involved. Specifically, 
we consider the three reconstruction problems, \ellone{}, \ellp{} and \tv{}: 
\begin{align*}
(\ellone{})&\qquad\qquad\qquad\argmin_\im ~\|\im\|_{1\phantom{\text{TV}}} \; \text{subject to} \quad \sysmat \im = \sino,\\
(\ellp)&\qquad\qquad\qquad\argmin_\im ~\|\im\|_{1\phantom{\text{TV}}} \; \text{subject to} \quad \sysmat \im = \sino, \quad \im \geq 0,\\
(\tv)&\qquad\qquad\qquad\argmin_\im ~\|\im\|_{\text{TV}\phantom{1}} \; \text{subject to} \quad \sysmat \im = \sino.
\end{align*}
The first two are standard 1-norm minimization (the latter with non-negativity constraint enforced) for reconstruction of images sparse in the image domain. The last is TV minimization for sparsity in the gradient domain. The TV semi-norm is defined as
\begin{align*}
 \|\im\|_\text{TV} = \sum_{j=1}^N \|D_j \im\|_2,
\end{align*}
where $D_j$ is a finite-difference approximation of the gradient at pixel $j$. In this work we use forward differences and Neumann boundary conditions. 

In the idealized setup we are interested in the central property of \emph{recoverability}: an image is said to be recoverable (from its ideal synthetic data) if it is the unique solution to the considered reconstruction problem. For example, we say that an image $\imorig$ is recoverable by \ellone{} from data $\sino = \sysmat \imorig$ if $\imorig$ is the unique \ellone{} solution.
The fundamental question we are interested in is:
\begin{center}
\emph{How few samples are enough for recovery of an image of a given sparsity by SR reconstruction?} 
\end{center}
In other words, we want to study recoverability as function of sparsity and sampling levels.
In the present work we will develop and apply a systematic analysis tool known as phase-diagram analysis from the field of compressed sensing (CS) for this purpose in the setting of CT.

\subsection{Compressed sensing} \label{subsec:cs}

The field of CS addresses precisely the question of how few samples one can acquire and still \emph{provably} recover the image.
In general, obviously, we need $N$ linearly independent samples of an image $\im \in \mathbb{R}^N$ to recover $\im$.
What CS says is that if the image $\im$ is sparse then by taking the right kind of samples we can recover $\im$ by SR from fewer than $N$ samples. Furthermore, the more sparse $\im$ is, the fewer samples will suffice.
CS was initiated with the works of Donoho \cite{Donoho2006} and Cand\`es et al. \cite{CandesTao2005decoding,candes2006robust}.
Before the advent of CS, SR reconstruction using the $1$-norm had been used heuristically for reduced sampling in CT \cite{delaney1998globally,li2002accurate},
but the works of Donoho and Cand\`es sparked renewed interest and a new focus on guarantees of accurate reconstruction.

An important quantity for CS guarantees is the restricted isometry property (RIP), which is defined as follows. A matrix $\sysmat$ is said to satisfy the RIP of order $s$ if there exists a constant $\delta_s \in (0,1)$ such that for all $s$-sparse signals $\im$ it holds that
\begin{align}
 (1-\delta_s)\|\im\|_2^2 \leq \|\sysmat \im\|_2^2 \leq (1+\delta_s)\|\im\|_2^2.
\end{align}
An example of a RIP-based CS guarantee is (see e.g. \cite{candes2008introduction}): 
If a matrix $\sysmat$ satisfies the RIP with $\delta_{2s} < \sqrt{2} - 1$, then an $s$-sparse image $\im$ will be recovered by \ellone{} from data $\sino = \sysmat \im$.

The problem is then to identify matrices satisfying this, and unfortunately computing RIP-constants is in general NP-hard \cite{Tillmann2014}. An important class of matrices that admit RIP-results are the Gaussian sensing matrices, for which matrix elements are independent samples from the zero-mean, unit-variance normal distribution. 
If the number of measurements $m$ satisfies
\begin{align}
 m \geq C \cdot s \cdot \log (N/s), \label{eq:ripnumsamples}
\end{align}
where $C$ is a constant, then with high probability a Gaussian sensing matrix possesses the RIP, such that the image $\im$ will be recovered. 

In a certain sense the Gaussian sensing matrices constitute an \emph{optimal sampling strategy} \cite{CandesTao2006nearoptimal,candes2008introduction}, because no other matrix type can provide the same recovery guarantee for fewer samples than 
\eqref{eq:ripnumsamples}. The importance of the Gaussian sensing matrices in CS is further established by many additional guarantees based for example on incoherence of the sensing matrix. It is not our intention to give a comprehensive review of CS theory here; such can be found in many places, for example the recent overview by Foucart and Rauhut \cite{FoucartRauhut:2013}. 

The prominent role of the Gaussian sensing matrices and other random matrix constructions in CS gives the impression that random sensing is a key CS feature and it is 
tacitly assumed in the imaging community that random sensing provides superior recoverability
performance to that of structured sampling.
This assumption has even lead researchers
to investigate hardware implementations of random sampling for CT \cite{Brady_CStomo:14}.
However, more recently novel CS guarantees have appeared for certain \emph{non-random} matrices \cite{GilbertIndyk2010}, which may be a step toward reduced focus on random sampling, although these matrices are also quite far from CT.

It is generally well-understood \cite{elad2010sparse,FoucartRauhut:2013} that current CS theory does not cover deterministic sampling setups in real-world applications. For CT in particular Petra and Schn\"orr \cite{Petra:2009,PetraSchnoerr2014} showed that CS guarantees are extremely poor. 
The main sensing problem of CT is its fundamental nature of sampling the object by line integrals. Each line integral only samples a small part of the object, thus leading to sparse, highly structured and coherent CT sampling matrices. In contrast CS sensing matrices, such as the Gaussian, are dense, have random elements and are incoherent, and hence fundamentally different.
In other words, there remains a large gap between the empirically observed effectiveness of SR in CT and the mathematical CS guarantees of accurate recovery typically involving random matrices.

\subsection{Own previous work and contribution of present work}

We have recently been interested in analyzing SR in CT from a CS perspective \cite{Joergensen_TMI:2013,Joergensen_eqconpap_v2_arxiv:2014,uniqueness_arxiv_2014}.
More specifically, 
we have studied recoverability from fan-beam CT data by 1-norm and TV regularization. We introduced the use of certain phase diagrams from CS to the setting of CT for systematically studying how recoverability depends on sparsity and sampling. Our work demonstrated quantitatively that recoverability from equi-angular fan-beam CT data for certain classes of test images exhibits a phase-transition phenomenon very similar to what is has been proved in CS for the Gaussian sensing matrices, as will be explained in Sec.~\ref{sec:phasediagramanalysis}.

In the present work we will further refine the phase-diagram analysis we introduced in \cite{Joergensen_eqconpap_v2_arxiv:2014,uniqueness_arxiv_2014} and demonstrate how it can be used to systematically provide quantitative insight of the undersampling potential of SR in CT by applying it to 3 cases. First, in Sec.~\ref{sec:phasediagramanalysis} we will give the sufficient theoretical background on phase-diagram analysis and the application to CT. Following that, we address in Sec.~\ref{sec:comparegaussian}, \ref{sec:randomsampling} and \ref{sec:large-scale} the following studies:
\begin{itemize}[leftmargin=0.7in]
 \item[Study A:] 
 How does CT-sampling compare in terms of recoverability to an optimal CS sampling strategy, i.e., using Gaussian sensing matrices? 
 
 \item[Study B:] Is recoverability improved by taking random CT measurements?
 
 \item[Study C:] How accurately can small-scale synthetic-data phase diagrams predict sufficient sampling for realistically-sized images of real objects? 
\end{itemize}
Finally in Sec.~\ref{sec:conclusion} we conclude the paper.

The purpose of Study A is to put the CT phase-transition behavior we observed in \cite{Joergensen_eqconpap_v2_arxiv:2014,uniqueness_arxiv_2014} more clearly into context of CS-theory. Quite surprisingly our results demonstrate that standard CT sampling is almost comparable with Gaussian sensing matrices in terms of recoverability. This is surprising since the Gaussian sensing matrices form an optimal CS sampling strategy, as explained previously in this section.

Study B addresses the use of random sampling in CT for potentially allowing for accurate reconstruction from fewer measurements than regular structured CT sampling. By use of phase-diagram analysis we will show that random sampling does \emph{not} lead to improved performance, but rather unchanged or in some cases even substantially reduced performance.

The purpose of Study C is to establish a connection to real-world CT image reconstruction by investigating the practical utility of phase diagrams for predicting how much CT data to acquire for reconstructing accurately a large-scale image of a given sparsity. 

In all three studies we use phase-diagram analysis as the main tool. Our goal is both to arrive at the particular insights of the three studies and to demonstrate phase-diagram analysis as a useful tool for systematically gaining quantitative understanding of SR in CT.

\section{Phase-diagram analysis} \label{sec:phasediagramanalysis}

\subsection{Theoretical phase-transition results}
As explained in Sec.~\ref{sec:introduction}\ref{subsec:cs} the Gaussian sensing matrices play a central role in CS. It is also possible to give a theoretical description of its \ellone{} and \ellp{} recoverability in terms of phase-diagram analysis. We present two different theoretical analyses, by Donoho and Tanner (DT) and by Amelunxen, Lotz, McCoy and Tropp (ALMT).

\begin{figure}[tb]
\centering
 \includegraphics[width=0.4\linewidth]{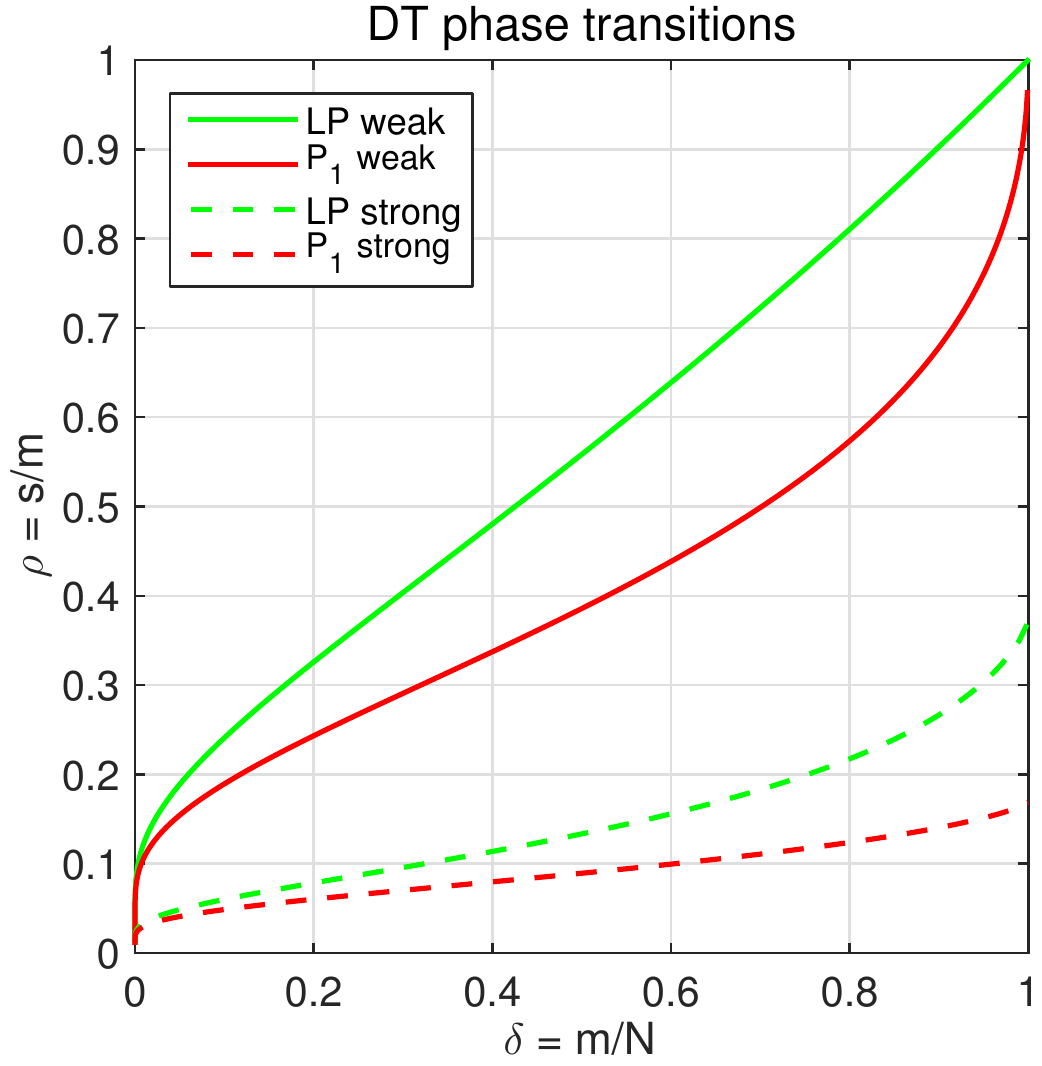}\imsep
 \includegraphics[width=0.4\linewidth]{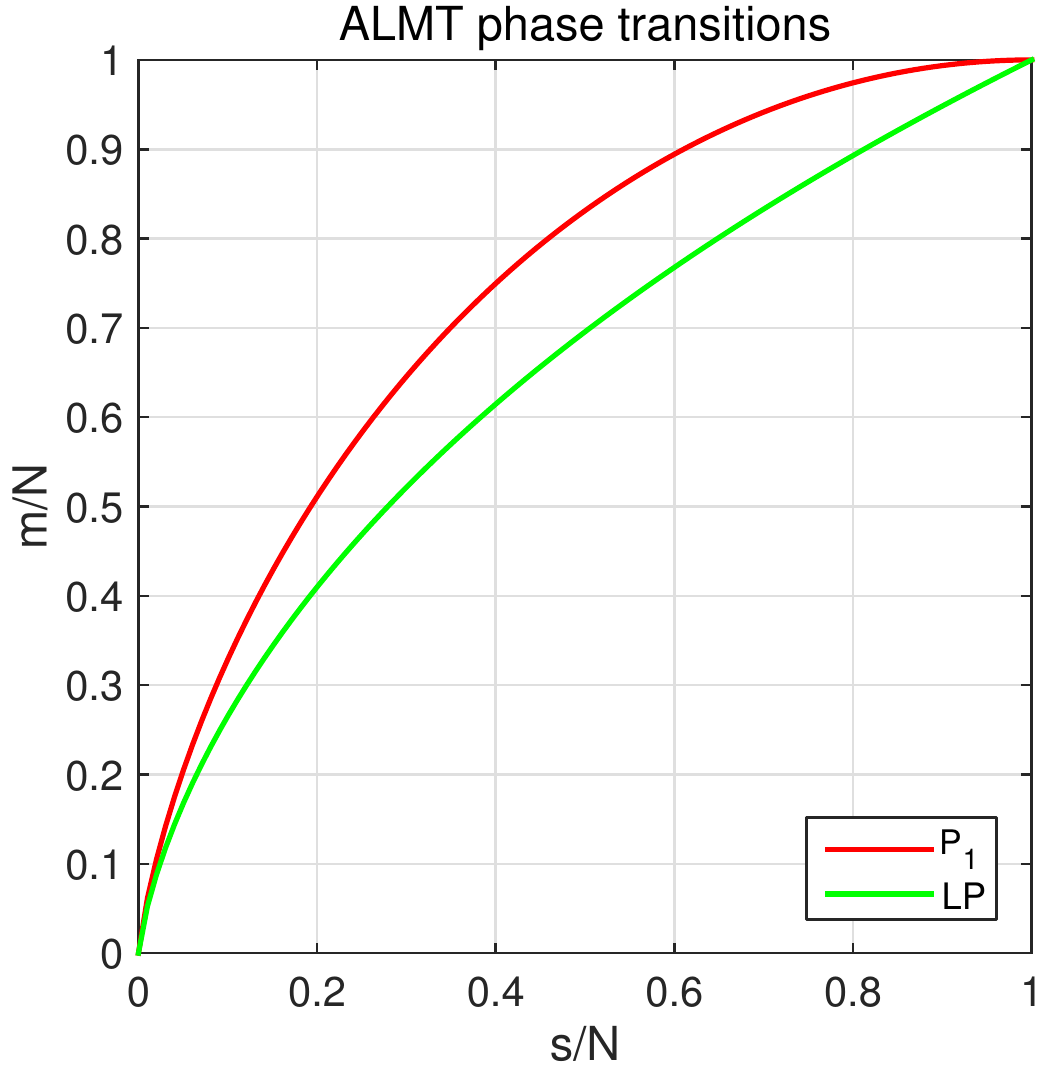}
 \caption{Theoretical phase-transition curves for Gaussian sensing matrices. Left: Donoho-Tanner (DT) asymptotic phase-transition curves for strong and weak recovery by \ellone{} and \ellp{}; recovery occurs \emph{below} the curves. Right: Amelunxen-Lotz-McCoy-Tropp (ALMT) phase-transition curves for recovery by \ellone{} and \ellp{}; recovery \emph{above} the curves.}\label{fig:DTALMTcurves}
\end{figure}

DT established in a series of papers \cite{DonohoTanner_nonneg:2005,DonohoTanner:2009,DonohoTanner_radically:2009,DonohoTanner_finitesize:2010} phase-transition behavior of the Gaussian sensing matrices. Their analysis is based on so-called neighborliness of random polytopes and builds on earlier work by Vershik and Sporyshev \cite{Vershik1992}. For an $s$-sparse signal $\im \in \mathbb{R}^N$ and $m$ samples, the DT phase diagram displays recoverability as function of $(\rho, \delta)$ for $\rho = s/m \in [0,1]$ and $\delta = m/N \in [0,1]$. For the set of $s$-sparse signals DT consider two notions of recoverability: strong, meaning that \emph{all} $s$-sparse signals are recovered, and weak, meaning that \emph{most} $s$-sparse signals are recovered at a given sampling level. DT then showed for the Gaussian sensing matrices and \ellone{} and \ellp{} that asymptotically 
there exist strong/weak phase-transition curves $\rho(\delta)$ such that at a sampling level of $\delta$ with 
high probability all/most 
signals with $\rho < \rho(\delta)$ 
will be recovered. Similarly,  with high probability all/most signals with $\rho > \rho(\delta)$ will fail to be recovered. The strong and weak phase-transition curves for \ellone{} and \ellp{} are shown in Fig.~\ref{fig:DTALMTcurves}~(left). Each curve partitions the phase space in two regions, one of full recovery (below the curve) and one of no recovery (above). We note that the weak full-recovery regions are substantially larger than their strong counterparts and that \ellp{} has a larger full-recovery region than \ellone{}. Both observations intuitively make sense. 
As we will demonstrate in Sec.~\ref{sec:comparegaussian}, the asymptotic weak phase-transition curves are in excellent agreement with empirical phase diagrams for finite-sized problems.

ALMT use a completely different analysis \cite{Amelunxen_arxiv:2014} based on the so-called statistical dimension of descent cones to prove non-asymptotic phase-transition behavior for the Gaussian sensing matrices. The ALMT phase diagram shows recoverability as function of $(s/N, m/N) \in [0,1]^2$. ALMT give phase-transition curves i.e. critical sampling values $m/N$ as function of sparsity values $s/N$ such that most images of a given sparsity are recovered from more samples than the critical level, and not recovered from fewer samples. The \ellone{} and \ellp{} ALMT phase-transition curves are shown in Fig.~\ref{fig:DTALMTcurves}~(right). Contrary to the DT phase-transition curves, the full recovery regions are above the curves. We will demonstrate in Sec.~\ref{sec:comparegaussian} that the ALMT phase-transition curves are in excellent agreement with empirical phase diagrams.

Regarding recovery guarantees for TV, we are only aware of the RIP-results by Needell and Ward \cite{needell2013cstv,NeedellWard:SIR:13}. 
To our knowledge it is an open question whether theoretical phase-transition results can be obtained. In the present work we demonstrate empirically that such behavior can be observed both from Gaussian and fan-beam CT sensing matrices.

In addition to the Gaussian sensing matrices, phase-transition behavior has been observed \cite{DonohoTanner:2009} for several other classes of random matrices and some theoretical analysis has been given \cite{Bayati:2012}. However, it remains open to establish phase-transition behavior for matrices occurring in practical imaging applications such as CT. Our motivation for the present work is precisely to establish that at least empirically it is possible to observe phase-transition behavior in CT.

\subsection{Experimental procedure of empirical phase-diagram analysis}\label{subsec:experimental}

Even though no theoretical phase-transition results exist for CT we can construct
empirical phase diagrams by repeatedly solving the same reconstruction problem over an ensemble of problem realizations for a range of sparsity and sampling levels. In our case we found that $100$ realizations at each sparsity and sampling level were enough to demonstrate phase-transition behavior. 

Each problem realization is generated in the following way: Given sparsity and sampling levels a test image $\imorig$ is generated, a sampling matrix $\sysmat$ is set up, and ideal data $\sino = \sysmat\imorig$ is computed. From the data $\sino$ the appropriate reconstruction problem is solved and the reconstruction is denoted $\imsol$. Recovery is declared if $\imsol$ is sufficiently close numerically to $\imorig$; here we test whether the relative 2-norm error $\|\imsol-\imorig\|_2 / \|\imorig\|_2 < \epsilon$, for some choice of threshold $\epsilon$. For \ellone{} and \ellp{} we found $\epsilon = 10^{-4}$ to be suitable, while for \tv{} we use $\epsilon = 10^{-3}$, as the conic optimization problem is more difficult to solve accurately.

As in \cite{Joergensen_eqconpap_v2_arxiv:2014,uniqueness_arxiv_2014} we use the commercial optimization software MOSEK \cite{MOSEK} to solve reconstruction problems required to construct a phase diagram. MOSEK uses a state-of-the-art primal-dual interior-point method, which allows us to solve \ellone{} and \ellp{} (recast as linear programs) and \tv{} (recast as a conic program), very accurately. An accurate solution is necessary for correctly assessing numerically whether an image is recovered, since numerical inaccuracies and approximate solutions may lead to the wrong decision. While allowing for high accuracy, interior-point methods are not efficient for large-scale problems.
For the reconstruction problems in Study C we use a large-scale optimization algorithm, which will be described there.

For the Gaussian sensing matrices, each problem realization contains a new realization of the sampling matrix, while in the fan-beam CT case a single matrix (at each sampling level) is used throughout. This is because
in CT we really are interested in the performance of a fixed matrix, which is specified by the physical scanner geometry.

For the ALMT phase diagrams we use $39$ relative sparsity levels $s/N = 0.025, 0.050, \dots, 0.975$ and $26$ sampling levels, namely from $1$ to $26$ equi-angular projection views. At $26$ views, the matrix has size $3338 \times{} 3228$ and is full rank, such that any image, independent of sparsity, will be recovered. For the DT phase diagram we use the same $26$ sampling levels in combination with $32$ sparsity levels (relative to the sampling level), i.e., $\rho = s/m = 1/32, 2/32,\dots,32/32$. 

With 100 realizations at each sparsity and sampling level, a total of 101,400 reconstruction problems need to be solved for a single ALMT phase diagram (at the chosen resolution), while the same number for a DT phase diagram is 83,200. Even with the small images used in this paper, our results have taken many hours of computing time on a cluster at DTU Computing Center.

%%%%%%%%%%%%%%%%%%%%%%%%%%%%%%%%%%%%%%%%%%%%%%%%%%%%%%%%%%%%%%%%%%%%%
\section{Study A: How does CT compare to CS?} \label{sec:comparegaussian} %\label{sec:smallresults}

As we have explained, the Gaussian sensing matrices are central to CS, since they admit strong theoretical results 
and are shown to form an optimal sampling strategy. In this study we use phase-diagram analysis to compare recoverability of fan-beam CT with the Gaussian sensing matrices. We will show that despite the lack of CS guarantees for fan-beam CT, we can empirically observe almost comparable recoverability.

\subsection{Measurement matrices}

We consider two types of measurement matrices: the Gaussian sensing matrices and a system matrix corresponding to a 2D equi-angular fan-beam scanning geometry. A Gaussian sensing matrix 
is generated by drawing
independent, identically distributed elements from the standard zero-mean unit-variance normal distribution. 

The 2D fan-beam CT system matrix is practically the same one we used in \cite{Joergensen_eqconpap_v2_arxiv:2014, uniqueness_arxiv_2014}, where it is described in detail, and the non-zero structure and the scanning geometry are illustrated in \cite{uniqueness_arxiv_2014}. In brief, we consider a disk-shaped image of $N$ pixels in total, inscribed in an $N_\text{side}\times{} \nside$ square pixel array. Fan-beam projections are recorded at $\nv$ equi-angular views of a $360^\circ$ scanning arc, each consisting of $2\nside$ pixels on a curved detector. The total number of measurements is $m=\nv\cdot 2 \nside$, and the $m\times{}N$ system matrix is computed by the function \texttt{fanbeamtomo} from the MATLAB$^\circledR$ toolbox AIR Tools \cite{Hansen2012}. The only difference from \cite{Joergensen_eqconpap_v2_arxiv:2014, uniqueness_arxiv_2014} is that the first angle is not chosen to be on a coordinate axes but offset by $20^\circ$. This offset regularizes the matrix by avoiding identical rows arising from rays 
in opposite views aligned with the coordinate axes.

\subsection{Image-domain sparsity}

\paragraph*{Signedspikes by \ellone{}}
We consider first the unconstrained problem \ellone{}. The standard image class considered in CS phase-diagram studies consists of images with random-valued pixels at random locations. We refer to this image class as signedspikes, see \cite{Joergensen_eqconpap_v2_arxiv:2014} for details and illustration. Specifically we generate a signedspikes image realization as follows: Given an image size (number of pixels) $N$ and sparsity (number of non-zero pixels) $s$, select uniformly at random $s$ pixels and assign values sampled from the uniform distribution on $[-1,1]$.

\begin{figure}[tb]
\centering
 \includegraphics[width=0.4\linewidth]{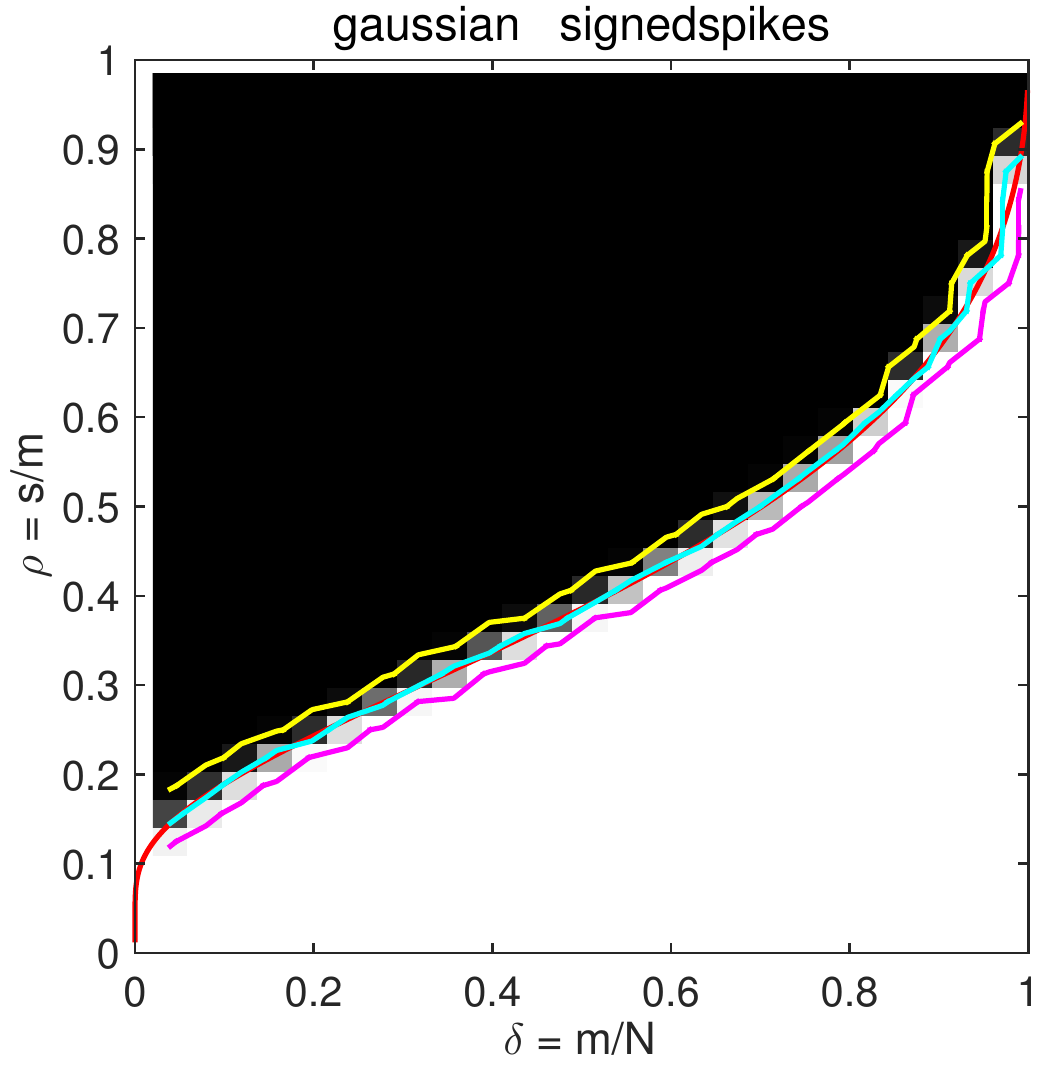}\imsep
  \includegraphics[width=0.4\linewidth]{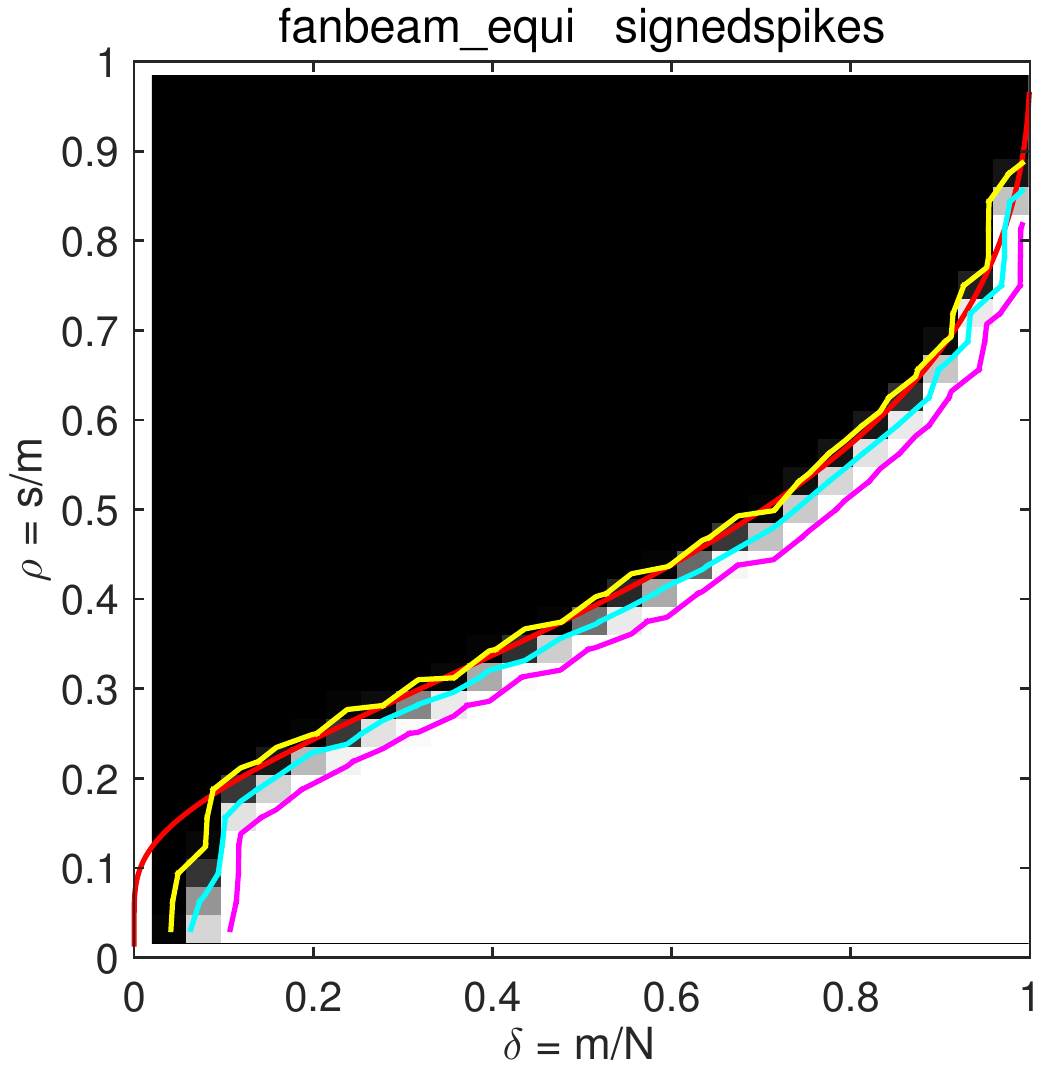}\\[0.1cm]
 \includegraphics[width=0.4\linewidth]{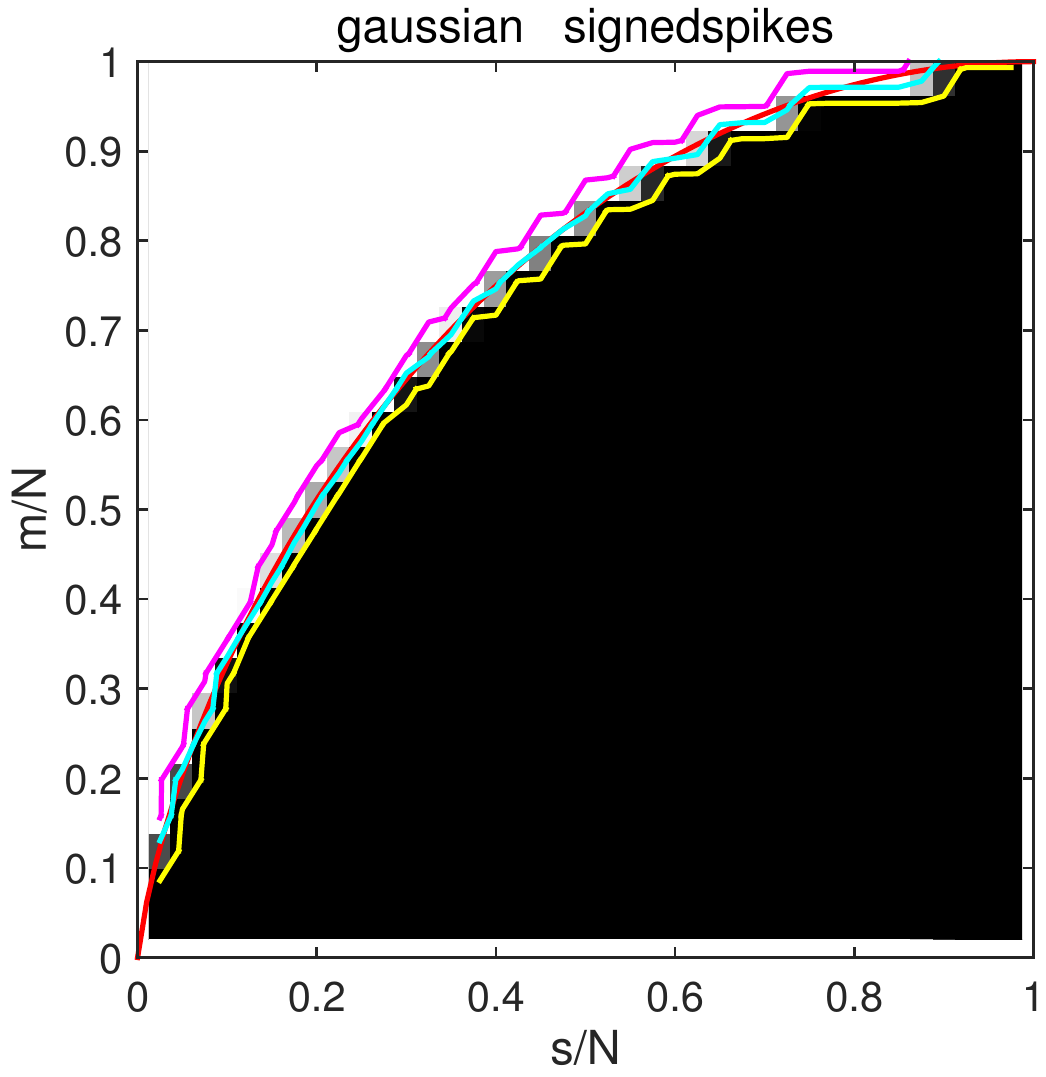}\imsep
  \includegraphics[width=0.4\linewidth]{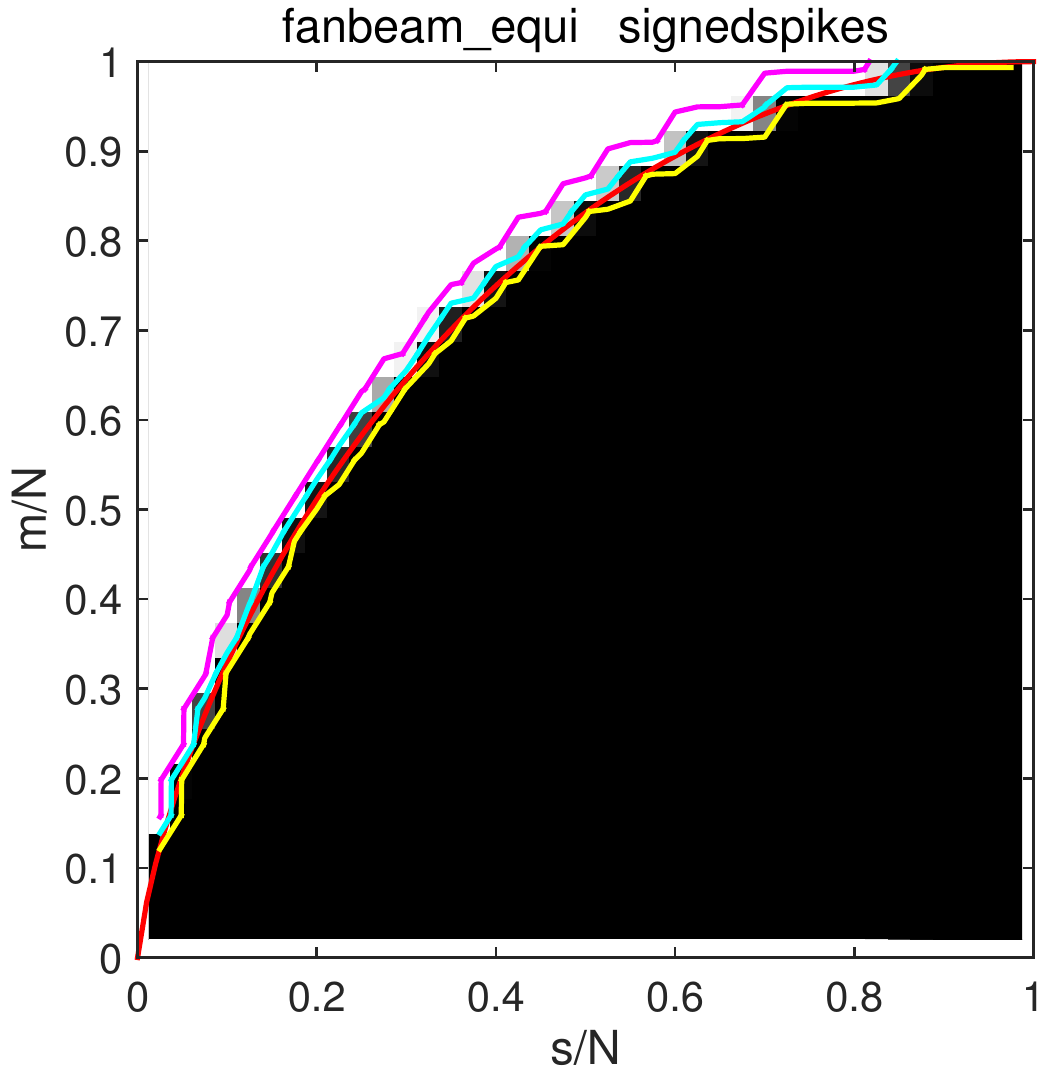}
   \caption{Phase diagrams for the signedspikes image class and \ellone{} reconstruction. DT phase diagrams (top row) and ALMT phase diagrams (bottom row). Gaussian sensing matrices (left) and fan-beam CT system matrices (right). Theoretical phase-transition curves for Gaussian sensing matrices (red), empirical phase-transition curve at $50\%$ contour line (cyan), $5\%$ and $95\%$ contour lines (yellow and magenta).\label{fig:dtalmt_signedspikes}}
\end{figure}

We generate DT and ALMT phase diagrams as described in Sec.~\ref{sec:phasediagramanalysis}\ref{subsec:experimental} for Gaussian and fan-beam CT sensing matrices, see Fig.~\ref{fig:dtalmt_signedspikes}.
At each sparsity and sampling level, the color represents the empirical success rate ranging from $0\%$ (shown black) to $100\%$ (shown white). Overlaid in cyan is the $50\%$ contour line indicating the empirical transition curve, as well as in yellow and magenta the $5\%$ and $95\%$ contour lines to quantify the transition width. Further, in red line is shown the theoretical phase-transition curve for the Gaussian sensing matrices.

We make the following observations. First, for the Gaussian sensing matrices, both the empirical DT and ALMT phase diagrams are in perfect agreement with the theoretical DT and ALMT phase-transition curves. This was to be expected but we include it here to verify that we can indeed reproduce the expected phase-transition curves using our software implementation. Second, and much more surprising, \emph{the fan-beam CT phase diagrams are almost identical to the Gaussian case.} The single apparent difference is in the bottom left corner of the DT phase diagram, where the CT recovery region does not extend to the same level as the Gaussian case. The poor CT recovery performance here is easily explained: the two leftmost columns correspond to a single projection and two projections $180^\circ$ apart, from which it is inherently difficult to produce an accurate reconstruction. Note that this issue is not apparent from the present ALMT phase diagram. Apart from this difference, the CT 
recovery performance is almost identical to the Gaussian case, in particular the transition is as sharp, as indicated by the $5\%$ and $95\%$ contour levels. On closer inspection the CT recovery region is slightly smaller than the Gaussian case, as seen by the lower cyan curve in the DT case and higher in the ALMT case. 

Nevertheless, considering that the Gaussian sensing matrices form an optimal sampling strategy, 
and that CT sampling matrices are highly structured, coherent and sparse, we find it extremely surprising to observe almost as good recoverability for CT.

\paragraph*{Non-negative spikes by \ellp{}}
\begin{figure}[tb]
\centering
 \includegraphics[width=0.4\linewidth]{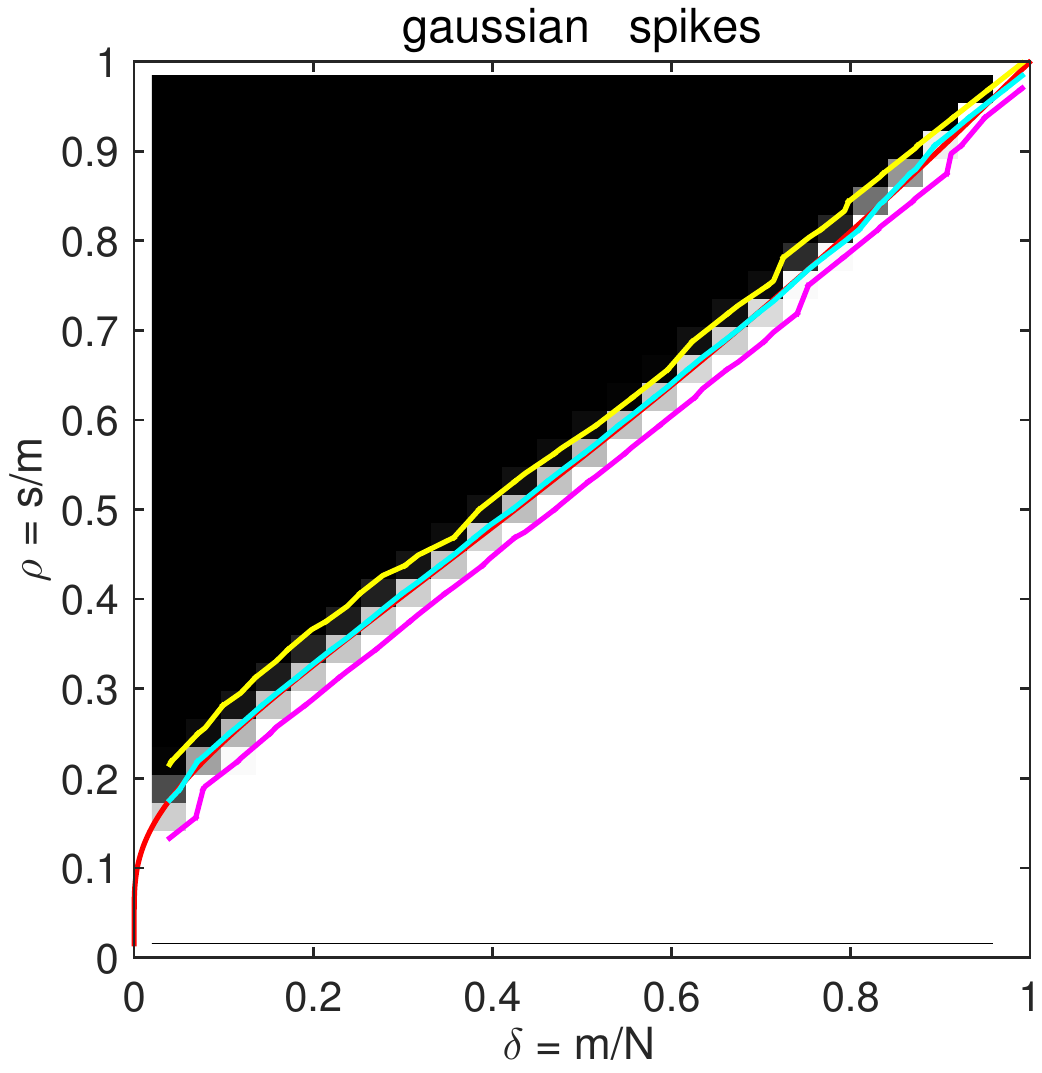}\imsep
  \includegraphics[width=0.4\linewidth]{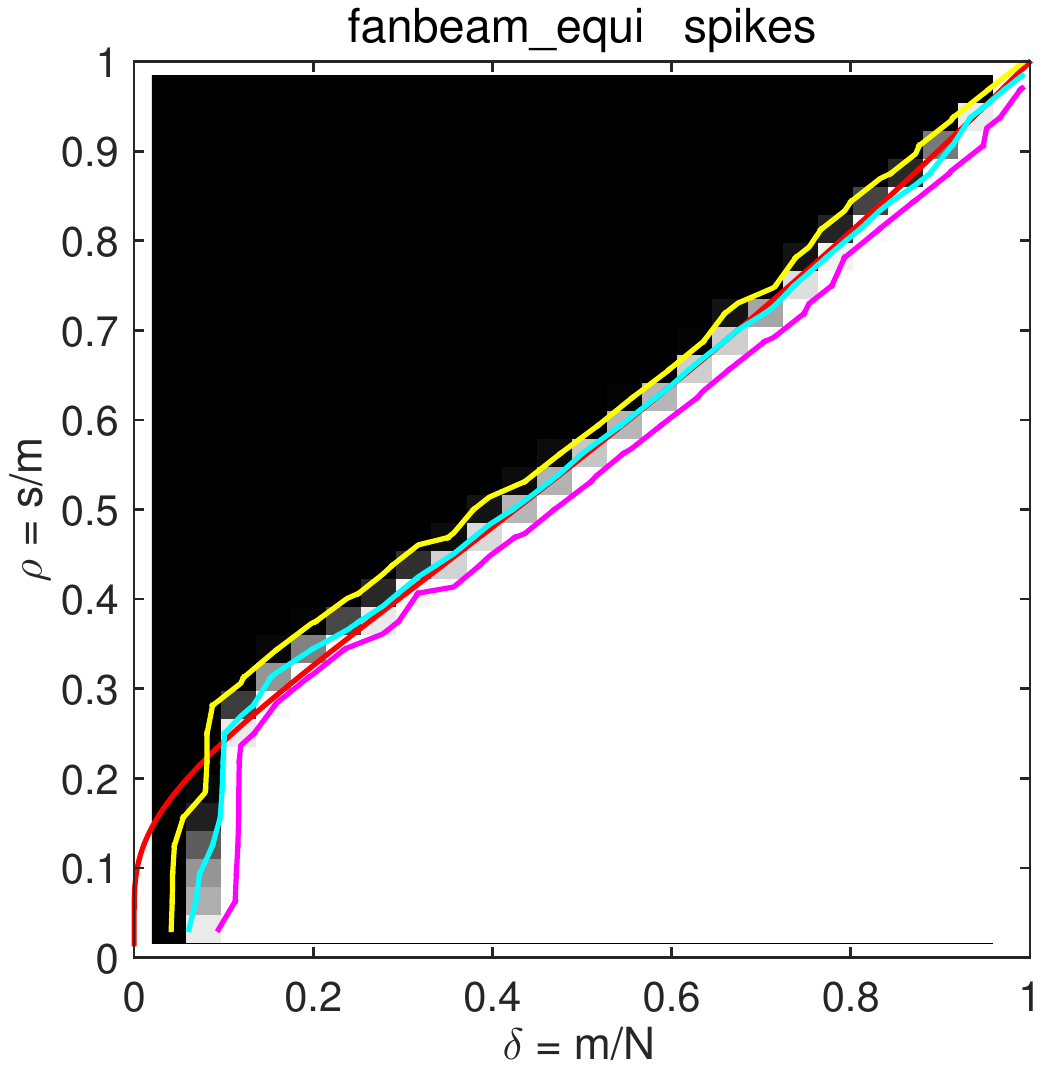}\\[0.1cm]
 \includegraphics[width=0.4\linewidth]{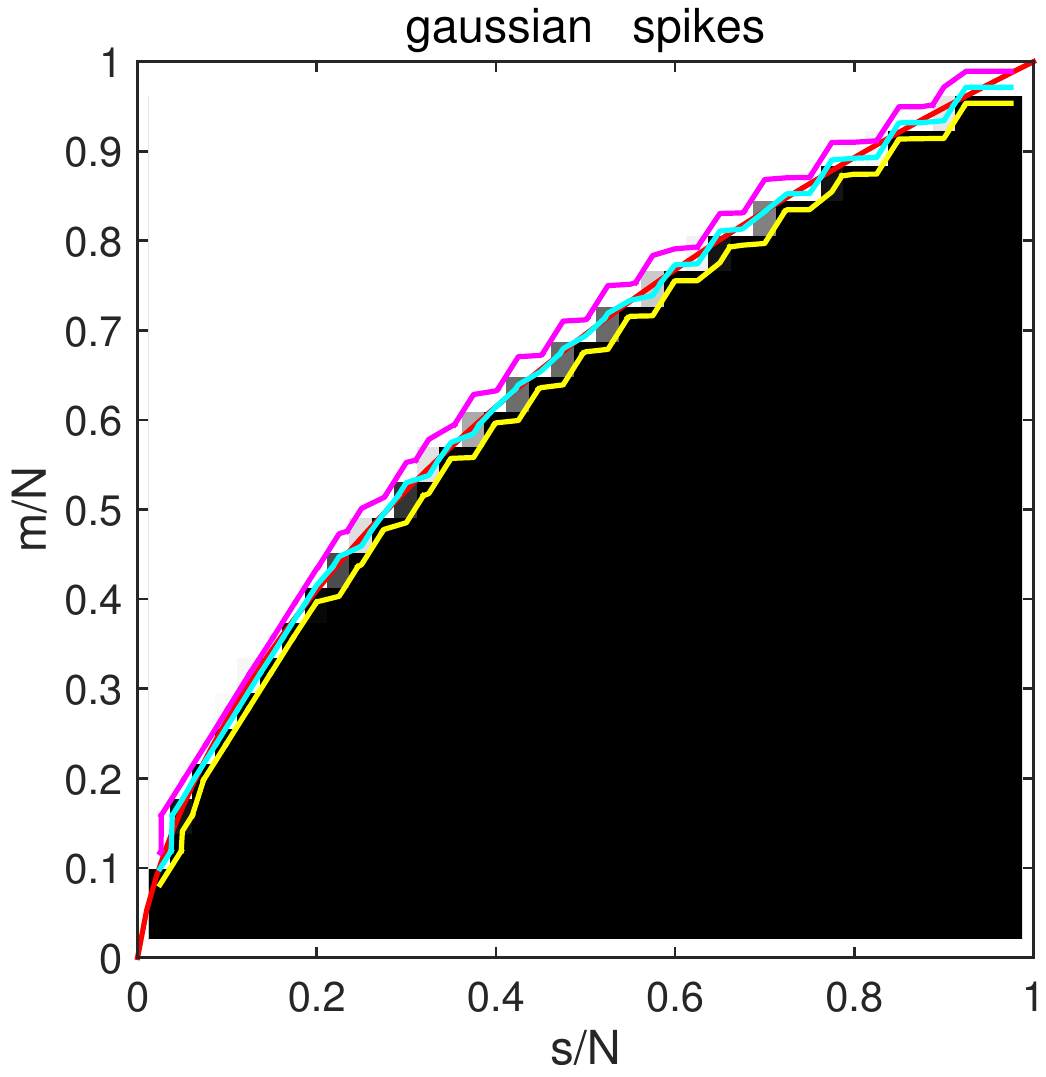}\imsep
  \includegraphics[width=0.4\linewidth]{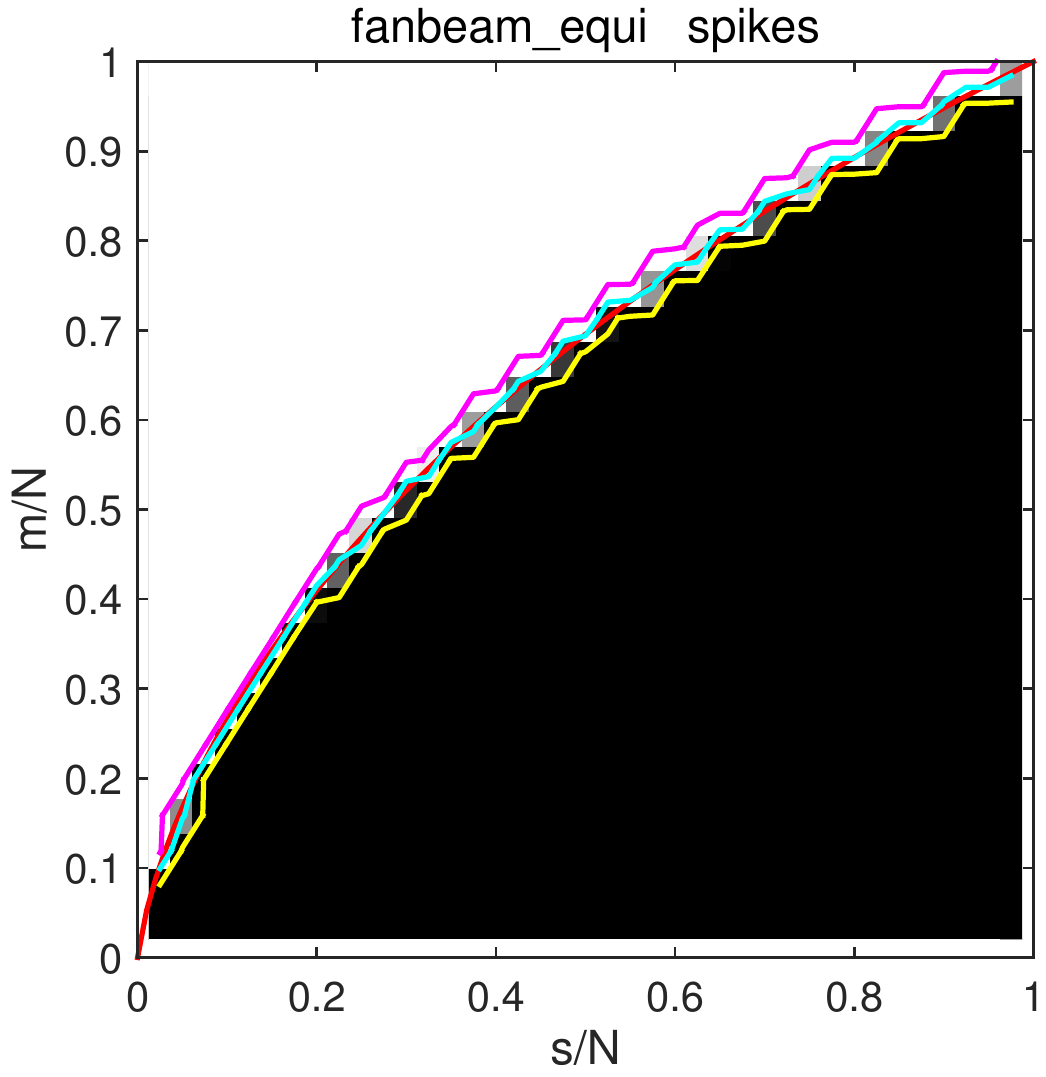}
   \caption{Phase diagrams for the non-negative spikes image class and \ellp{} reconstruction. DT phase diagrams (top row) and ALMT phase diagrams (bottom row). Gaussian sensing matrices (left) and fan-beam CT system matrices (right). Theoretical phase-transition curves for Gaussian sensing matrices (red), empirical phase-transition curve at $50\%$ contour line (cyan), $5\%$ and $95\%$ contour lines (yellow and magenta)\label{fig:dtalmt_spikes}}
\end{figure}
Typically in CT a non-negativity constraint can be employed since the imaged quantity, the linear attenuation coefficient, is non-negative, and hence the reconstruction problem \ellp{} is appropriate. 
For $\ellp$ we consider the natural non-negative version of the signedspikes class, which we call spikes, with the single change that 
values are sampled from the uniform distribution on $[0,1]$, see \cite{Joergensen_eqconpap_v2_arxiv:2014} for illustration.

We construct again empirical DT and ALMT phase diagrams and display them in Fig.~\ref{fig:dtalmt_spikes}
together with the theoretical Gaussian-case phase-transition curves for \ellp{}. Also in this case, the CT phase diagrams are almost identical to the Gaussian case, in terms both of the empirical phase-transition curve and the width as indicated by the $5\%$ and $95\%$ contour lines. In fact, the similarity is even larger as the cyan $50\%$ contour in the CT case coincides with the theoretical transition curve, except at the bottom-left corner of the DT phase diagram, as before caused by having only $1$ or $2$ CT projections.
In accordance with the theoretical curves, we see that even fewer samples suffice for recovery in the non-negative case compared to before.

\paragraph*{A structured image class}
CS recovery guarantees for example for the Gaussian sensing matrices state that the sufficient number of samples depends on the signal only in terms of the signal sparsity. That is, signals with structure in the non-zero locations should not require a different number of samples for recovery than unstructured signals such as the spikes images.  Does the same hold for for CT? We will demonstrate that the answer is no. Due to non-zero pixels selected at random in the spikes classes there is no structure, i.e., correlation between neighboring pixels. As an example of a class of sparse images with some structure in the non-zero locations we use the 2-power class from \cite{Joergensen_eqconpap_v2_arxiv:2014}. This image class is based on a breast tissue model, but for our purpose here, it suffices to say some correlation has been introduced between neighbor pixel values. 

\begin{figure}[tb]
\centering
 \includegraphics[width=0.4\linewidth]{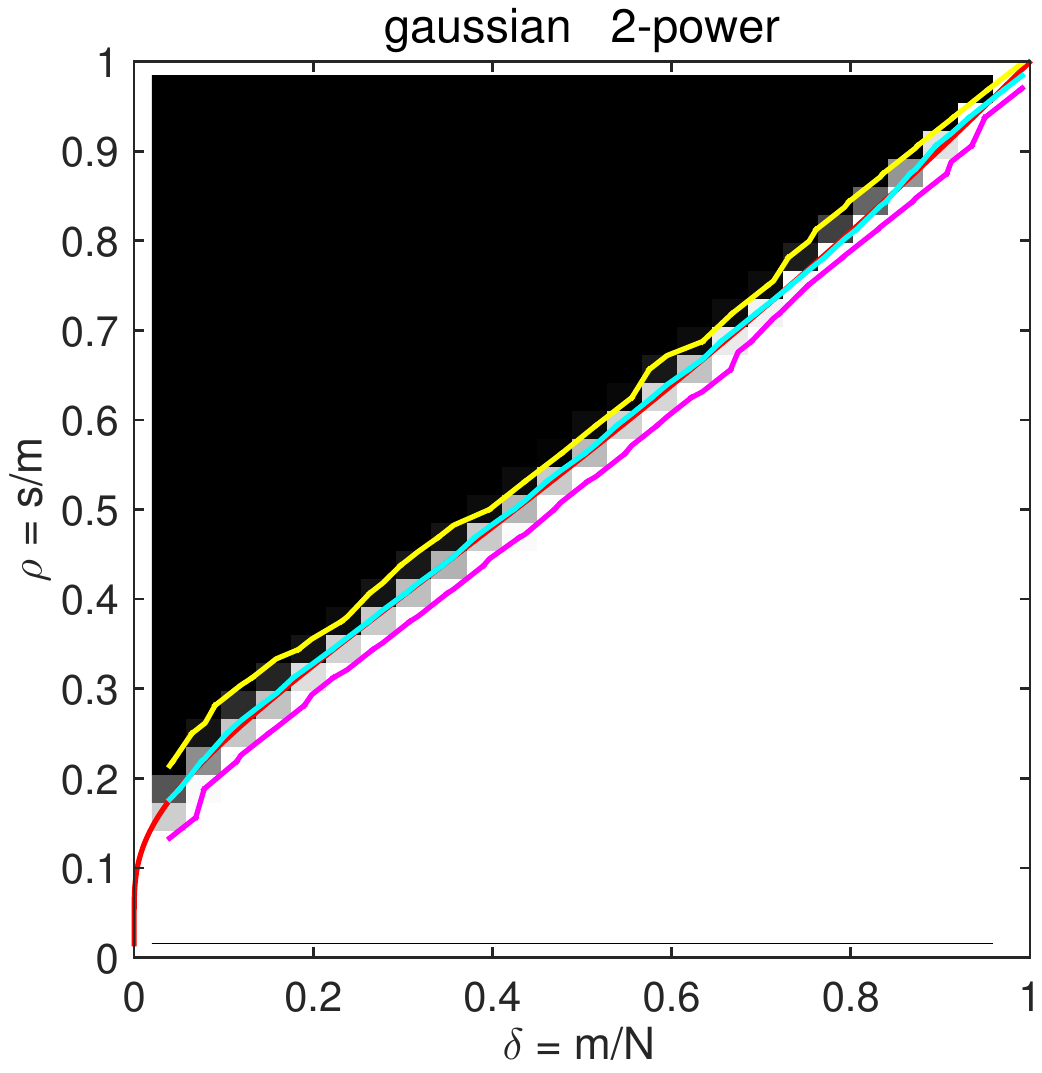}\imsep
  \includegraphics[width=0.4\linewidth]{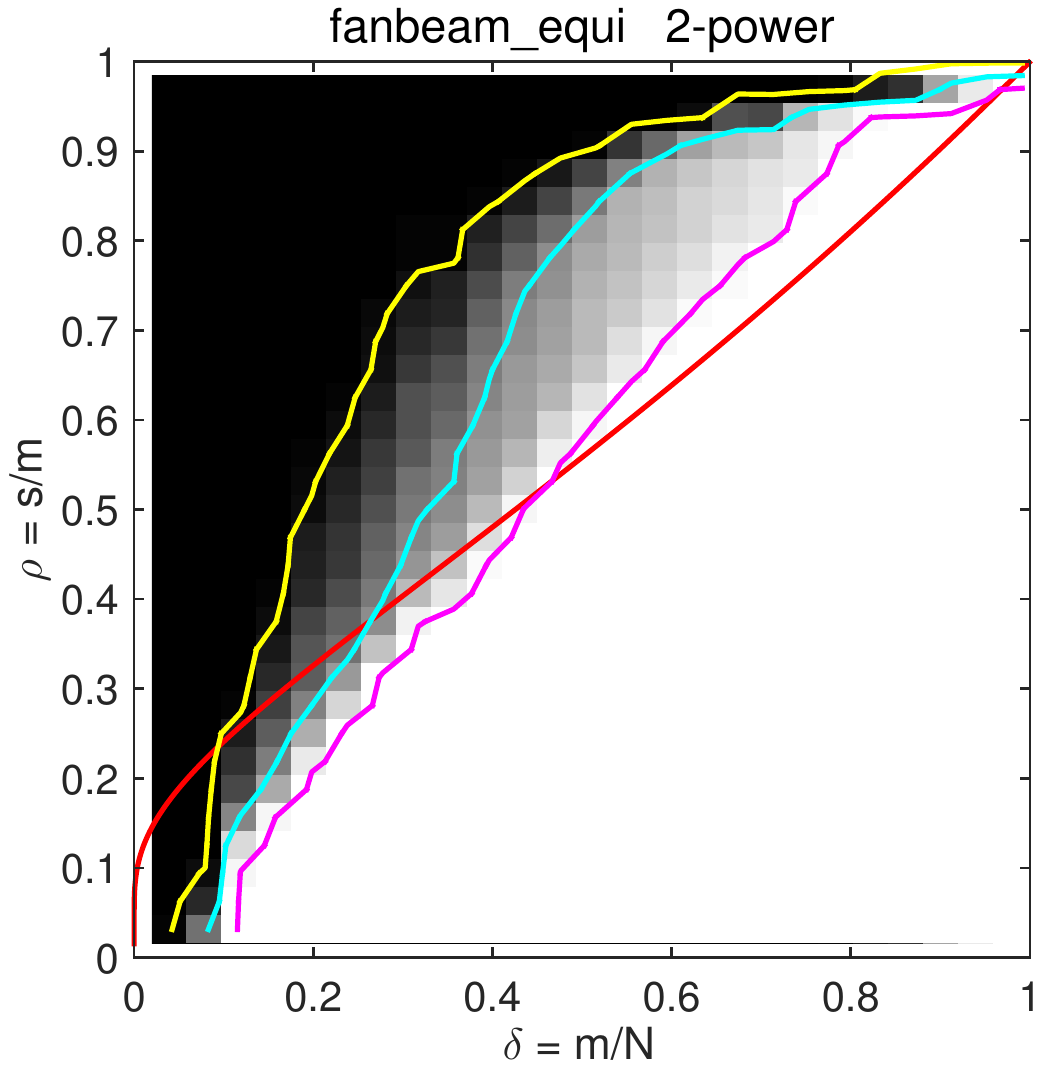}
   \caption{DT phase diagrams for the 2-power image class and \ellp{} reconstruction. Gaussian sensing matrices (left) and fan-beam CT system matrices (right). Theoretical phase-transition curves for Gaussian sensing matrices (red), empirical phase-transition curve at $50\%$ contour line (cyan), $5\%$ and $95\%$ contour lines (yellow and magenta).\label{fig:dt_fftpower_2_0}}
\end{figure}

Images from the 2-power class are non-negative, so we use \ellp{} for reconstruction, create DT phase diagrams, see Fig.~\ref{fig:dt_fftpower_2_0}, and compare with the spikes-class DT phase diagrams in Fig.~\ref{fig:dtalmt_spikes}, omitting ALMT phase diagrams for brevity. As expected, our results verify that image structure does not matter for the Gaussian sensing matrices, as the DT phase diagram is identical to the spikes case. But, for the fan-beam CT case the phase diagram has changed drastically, most notably the transition is now much smoother as indicated by the $5\%$ and $95\%$ contour lines. Also the empirical phase-transition curve ($50\%$ contour line) has moved away from the theoretical curve. We note that at low sampling (left part) the transition is lower while at high sampling, it is higher, so recoverability can be both better and worse, depending on sampling level.
The $95\%$ contour line limits a region of almost full recovery, and 
this region is not much different from the spikes case.

The 2-power result for CT is in stark contrast to the Gaussian sensing matrix behavior in Fig.~\ref{fig:dtalmt_spikes}.
We conclude that even though the spikes results suggest close resemblance of CT with the optimal CS case of Gaussian sensing matrices, the 2-power result makes it clear that CT is more complex.

\subsection{Gradient-domain sparsity}

Sparsity in the image domain is interesting due to well-developed theory, in particular for Gaussian sensing matrices. For CT it is more common to expect sparsity in the gradient domain, which has motivated the successful use of TV-regularization.
However, to the best of our knowledge, no phase-transition behavior has been proved, not even for the Gaussian case. Here, we demonstrate empirically that for both Gaussian and CT sensing matrices similar sharp phase transitions can be observed.

For generating images sparse in the gradient domain we use the image class from \cite{uniqueness_arxiv_2014} alternating-projection for (isotropic) TV, which we here refer to as altprojisotv. An image is generated in an iterative procedure of taking alternating projections onto the range of the gradient operator and thresholding the number of non-zeros in the image gradient to the desired sparsity level; see \cite{uniqueness_arxiv_2014} for details and illustration.

\begin{figure}[tb]
\centering
 \includegraphics[width=0.4\linewidth]{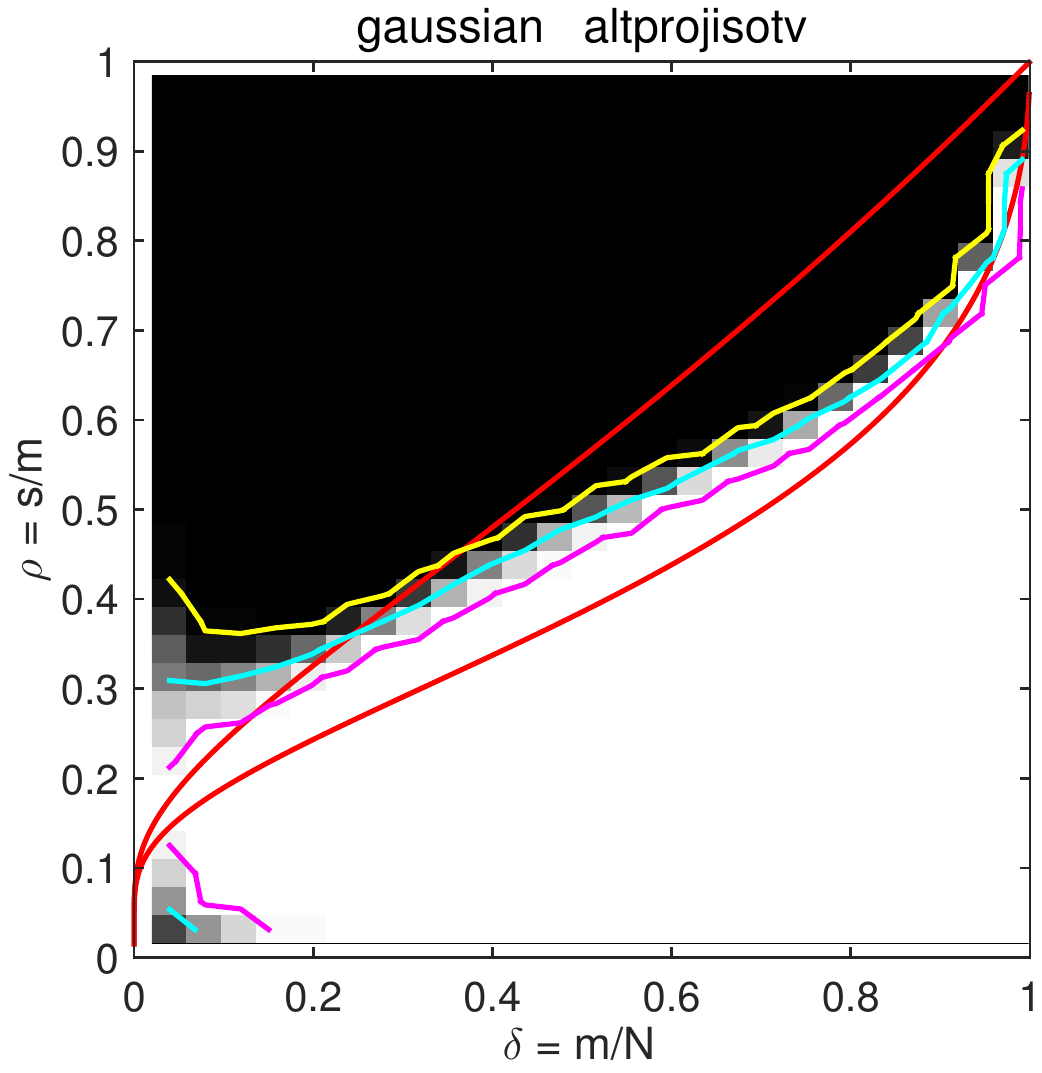}\imsep
  \includegraphics[width=0.4\linewidth]{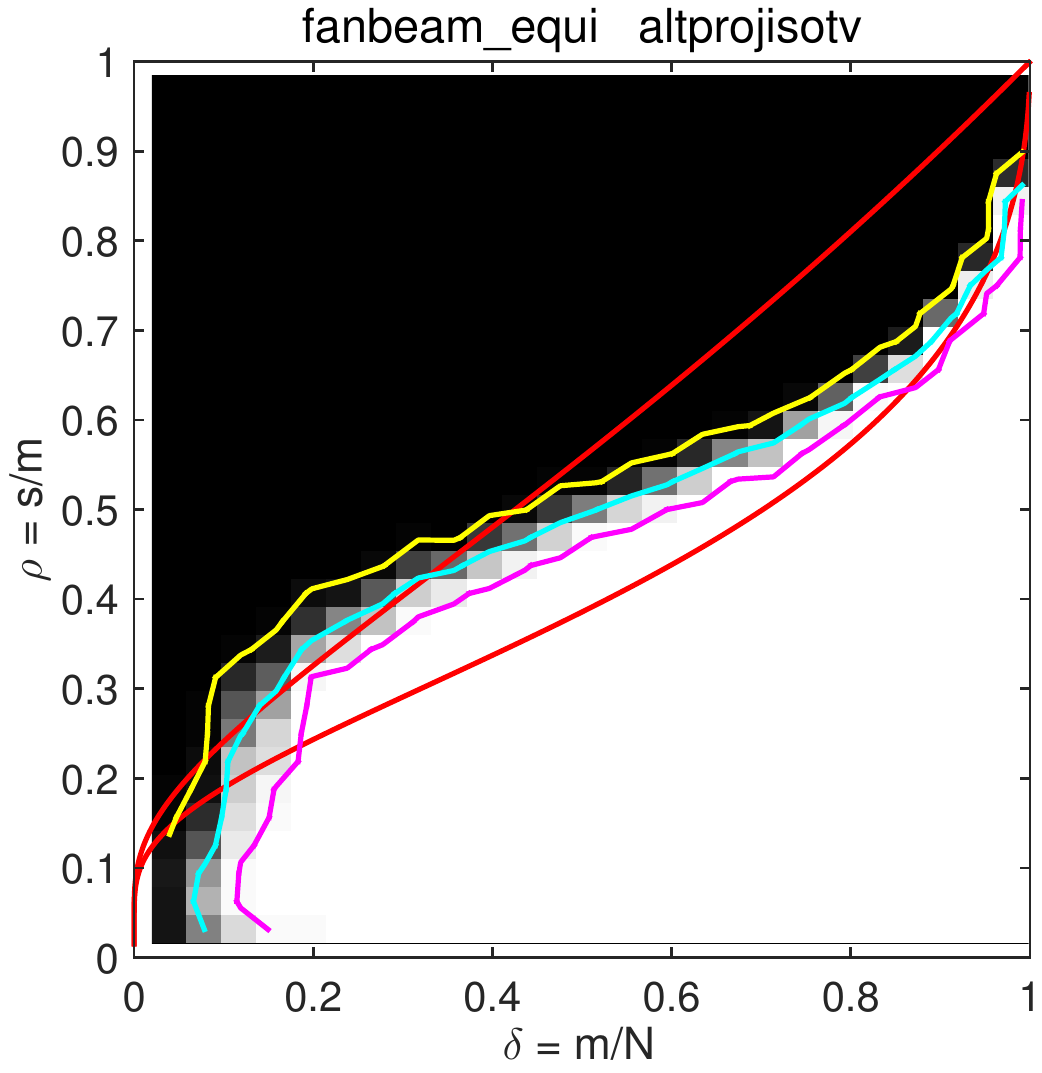}\\[0.1cm]
 \includegraphics[width=0.4\linewidth]{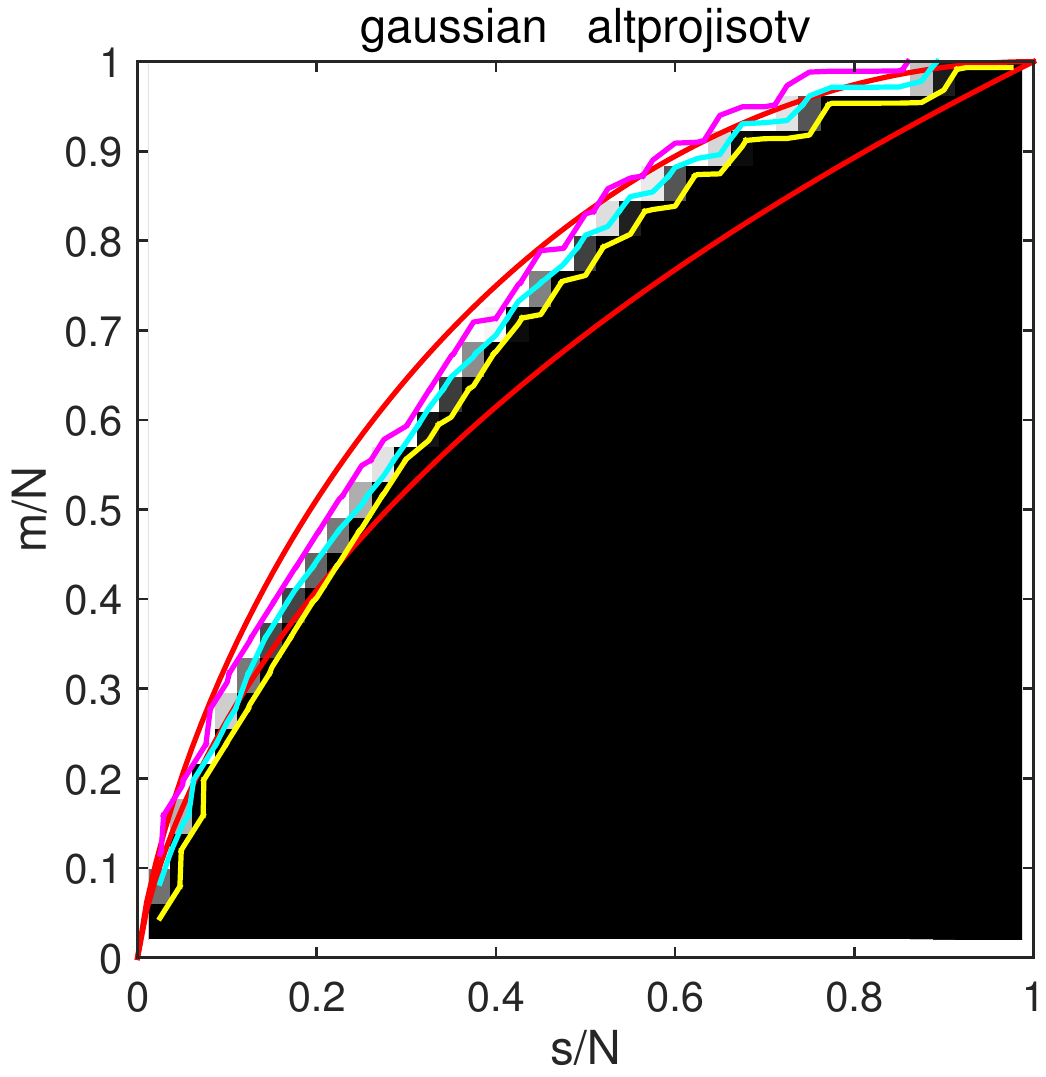}\imsep
  \includegraphics[width=0.4\linewidth]{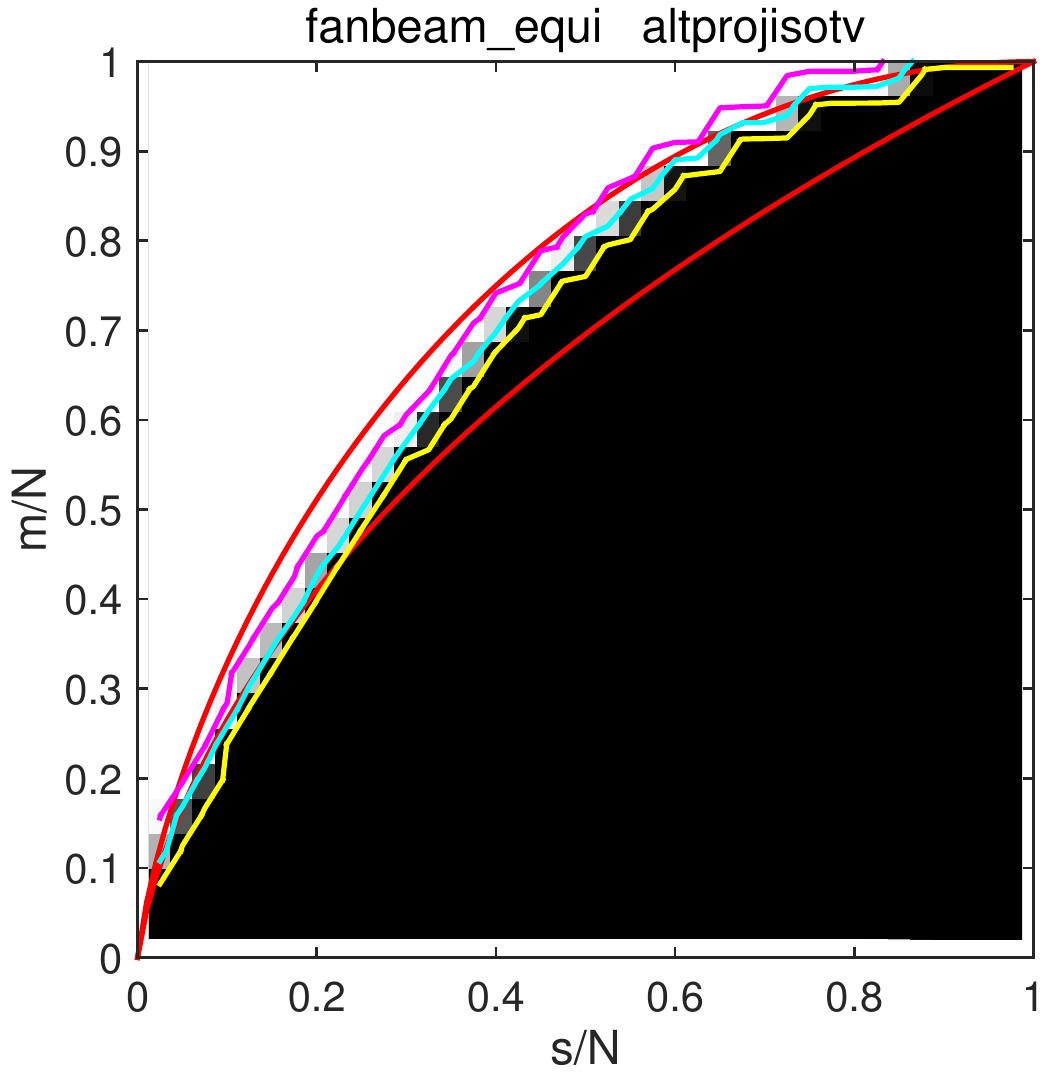}
      \caption{Phase diagrams for the altprojisotv image class and \tv{} reconstruction. DT phase diagrams (top row) and ALMT phase diagrams (bottom row). Gaussian sensing matrices (left) and fan-beam CT system matrices (right). Theoretical phase-transition curves for \ellone{} and \ellp{} reconstruction for Gaussian sensing matrices (red), empirical phase-transition curve at $50\%$ contour line (cyan), $5\%$ and $95\%$ contour lines (yellow and magenta).
      \label{fig:dtalmt_altprojisotv}}
\end{figure}

Once again, we construct DT and ALMT phase diagrams, see Fig.~\ref{fig:dtalmt_altprojisotv}; this time with sparsity values referring to gradient-domain sparsity. We observe also in this case a sharp phase transition both in the DT and ALMT phase diagrams. In the lack of a theoretical reference curve for TV we compare with the \ellone{} and \ellp{} curves and find that transition takes place between the two curves. 

An irregularity is observed in the bottom-left corner of both DT phase diagrams. The explanation is that 
the altprojisotv procedure has difficulty in generating images which are extremely sparse in the gradient domain.
In spite of the irregularity, we find that our empirical TV results convincingly demonstrate a that sharp phase transition takes place also in the TV case, dividing the phase space into regimes of full and no recovery, and again that CT recoverability is similar to the Gaussian case.

\subsection{Conclusion on Study A}
We used phase-diagram analysis to compare fan-beam CT recoverability with optimal CS-sampling using the Gaussian sensing matrices. For unstructured signed images with \ellone{} and non-negative images with \ellp{} we found almost identical phase-transition behavior in terms of critical sampling level and width of the transition. We thereby demonstrated that empirically fan-beam CT in the average-case performs close to the optimal. While recoverability by the Gaussian sensing matrices was unaffected by the introduction of structure in the non-zero pixels, fan-beam CT recoverability drastically changed to a much smoother transition. Interestingly, except for the lowest-sampling range, the recovery region actually became larger, meaning that many images at a given sparsity level are recovered from fewer samples than the Gaussian sensing matrices' critical sampling level. In spite of the close resemblance on the unstructured images, this example demonstrates that fan-beam CT is fundamentally different 
from 
the Gaussian sensing matrices.

Also in case of \tv{} recoverability we found almost identical behavior of fan-beam CT and the Gaussian sensing matrices. In particular in both cases we saw a sharp phase transition, thus suggesting that the phase-transition phenomenon generalizes to \tv{}. To our knowledge no theoretical explanation of this observation has been given in the literature.

%%%%%%%%%%%%%%%%%%%%%%%%%%%%%%%%%%%%%%%%%%%%%%%%%%%%%%%%%%%%%%%%%%%%%
\section{Study B: Is random sampling beneficial in CT?} \label{sec:randomsampling}

As mentioned in the introduction random sampling is an optimal strategy and important in many recovery guarantees. Sampling in CT is normally done in very structured manner and a natural contemplation is therefore whether the introduction of some form of randomness could lead to recovery guarantees for CT or improved recoverability compared to regular sampling.
In this study we use phase-diagram analysis to investigate whether CT sampling strategies involving randomness can improve the recoverability of sparse images, i.e., enable accurate reconstruction of images of a given sparsity from fewer measurements than regular equi-angular fan-beam CT.

\subsection{Measurement matrices}
Many forms of randomness can be conceived in CT sampling. In this work we consider two straightforward ones. First, a fan-beam geometry denoted fanbeam\_rand in which the source angular positions are no longer equi-distant but sampled uniformly from $[0,360^\circ]$. 
Second, we consider a setup we denote random\_rays of independent random rays through the image. Each ray is specified by two parameters: the angle of the ray with a fixed coordinate axis and the intersection of the ray with the orthogonal diameter of the disk-shaped image. The angle and intersection are sampled from the uniform distributions on $[0,180]^\circ$ and $[-N_\text{side}/2, N_\text{side}/2]$, respectively, where $N_\text{side}$ is the diameter length and image is assumed centered around the origin.

\subsection{Image-domain sparsity}
\begin{figure}[tb]
\centering
\newcommand{\ww}{0.4\linewidth}
 \includegraphics[width=\ww]{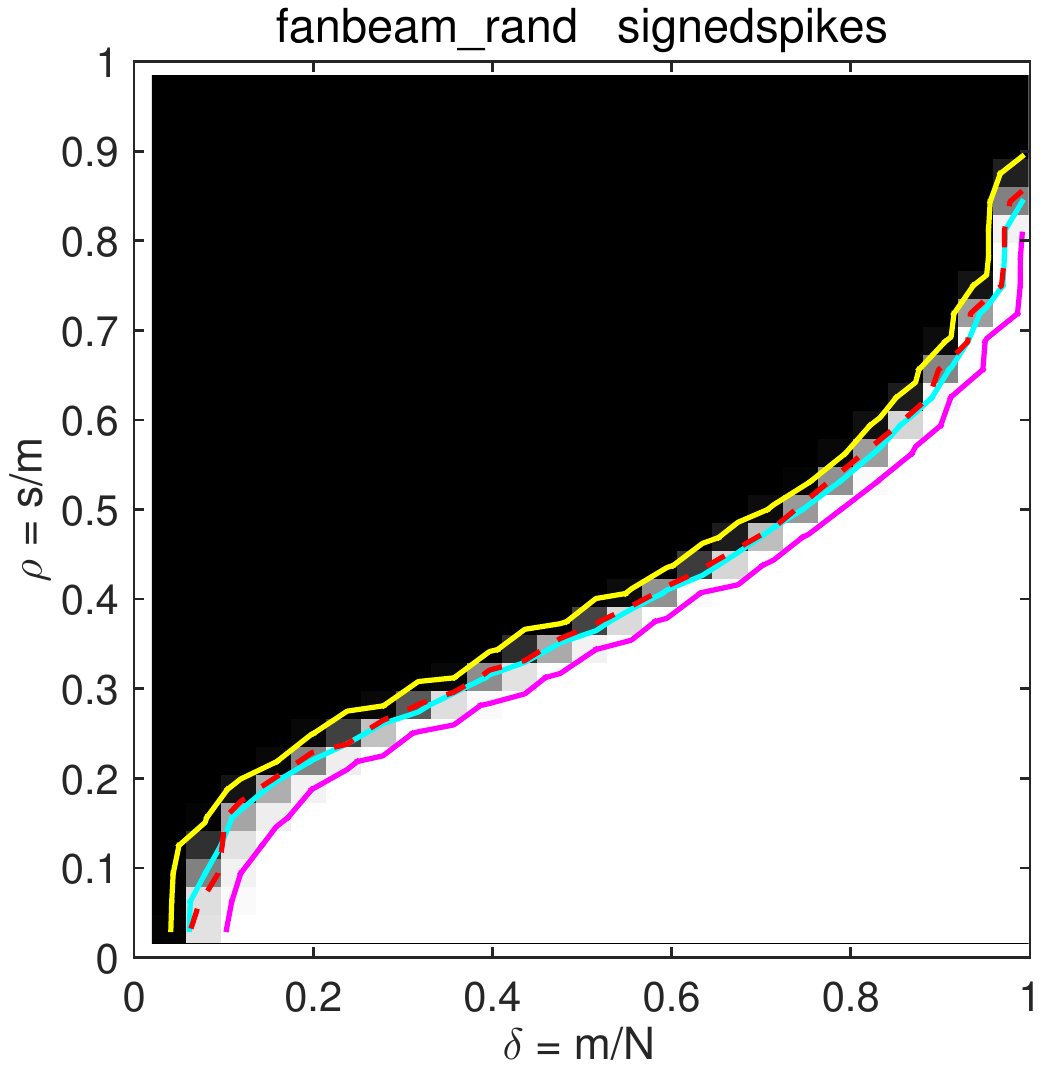}\imsep
  \includegraphics[width=\ww]{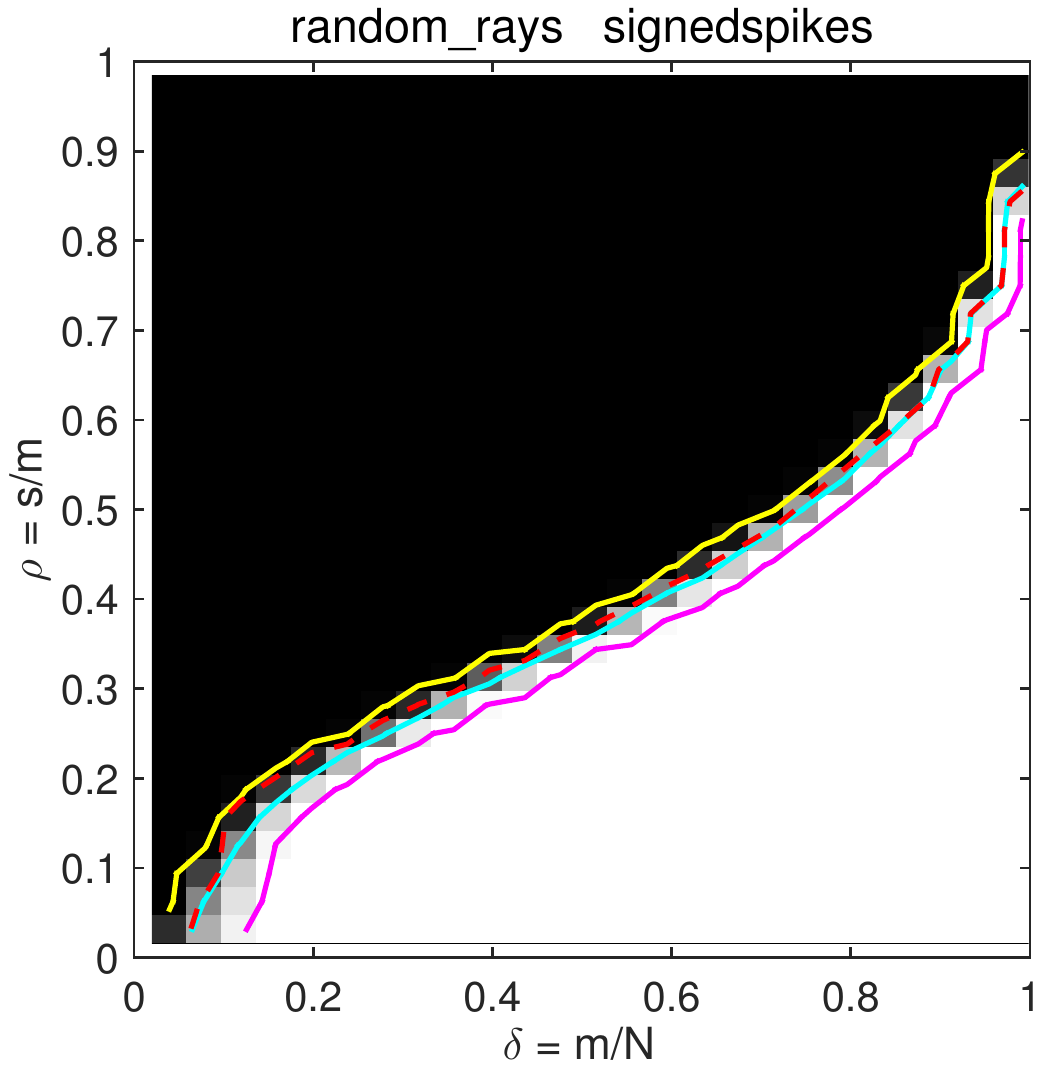}\\[0.1cm]
   \includegraphics[width=\ww]{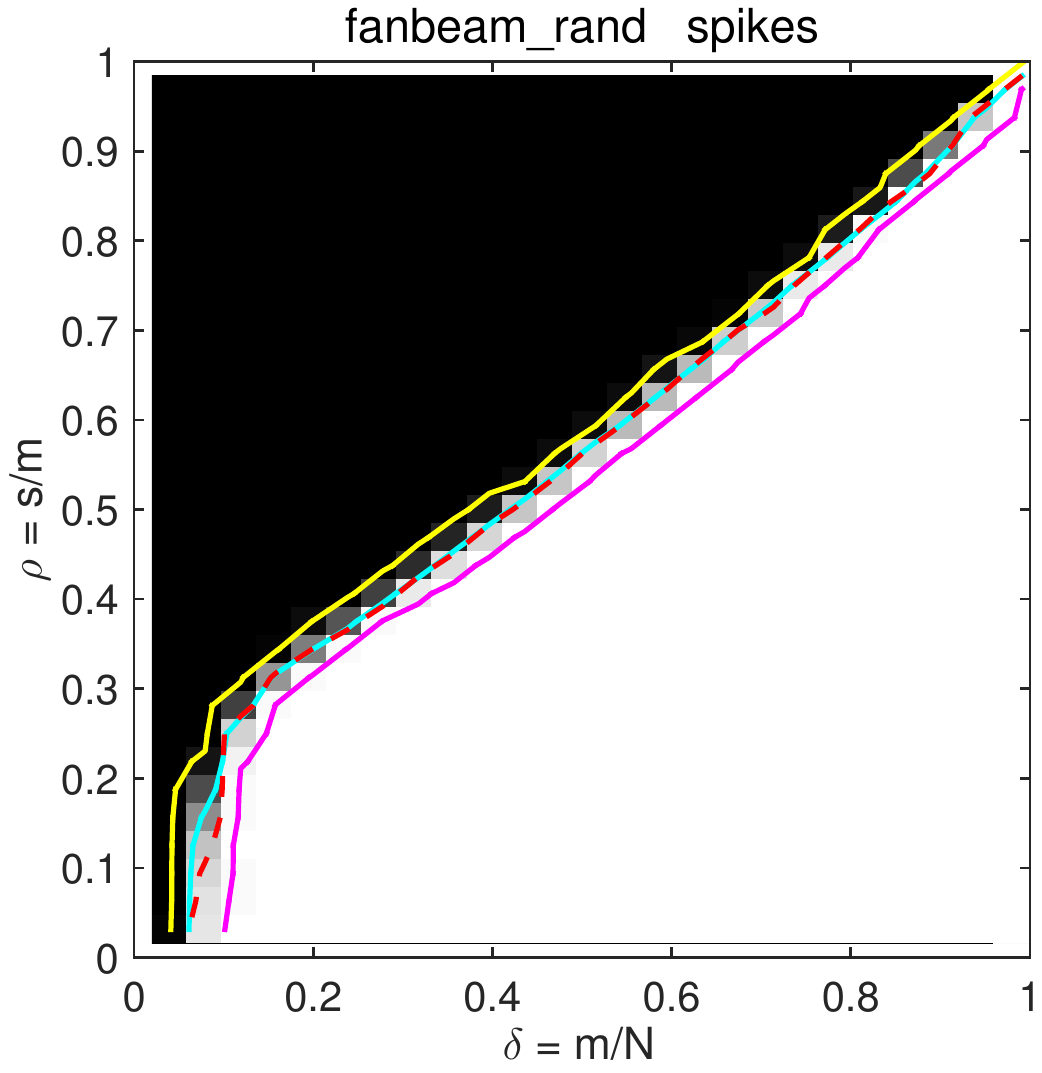}\imsep
  \includegraphics[width=\ww]{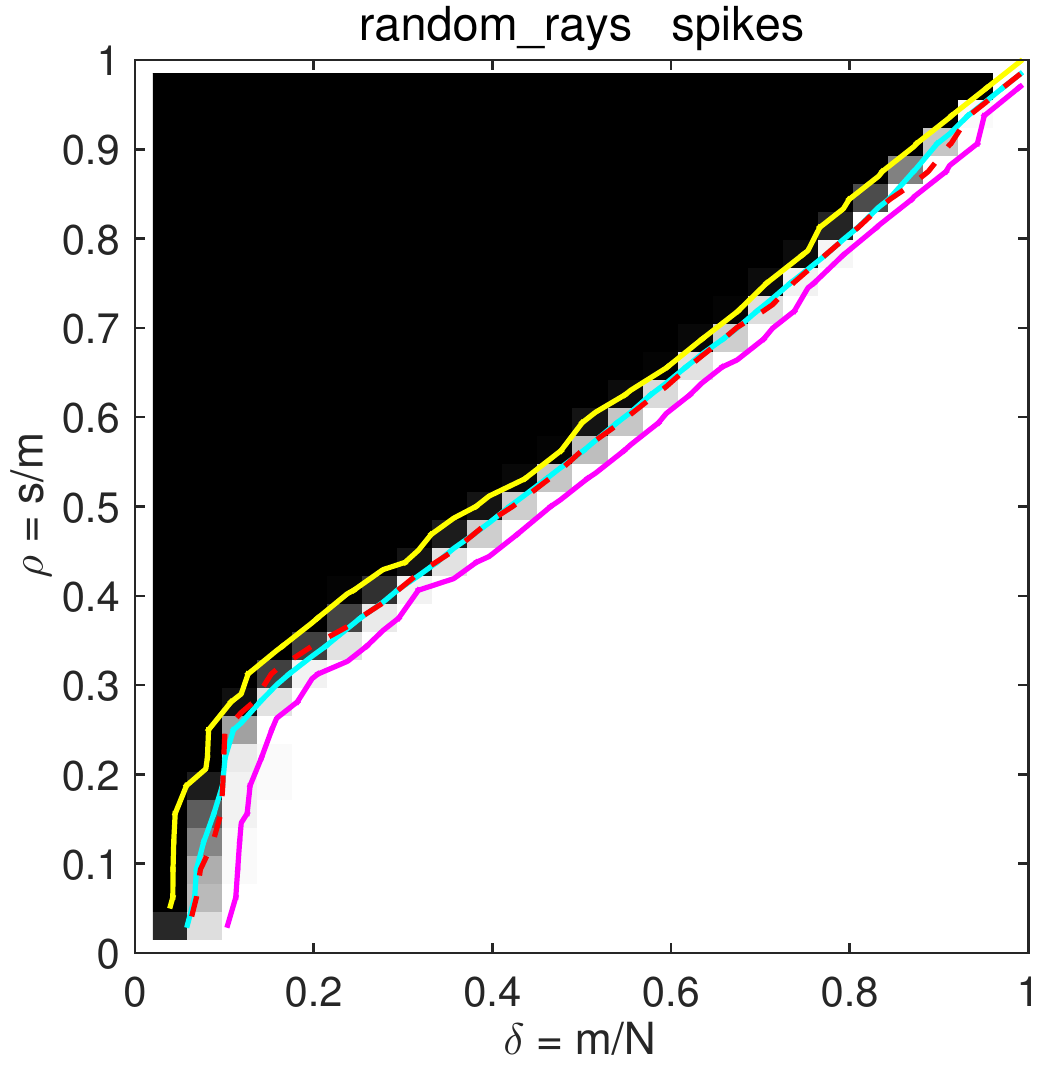}\\
\caption{DT phase diagrams. Signedspikes image class and \ellone{} reconstruction (top row) and spikes image class and \ellp{} reconstruction (bottom row). Fan-beam with random source positions (left) and random rays geometry (right). Empirical phase-transition curve for equi-angular fan-beam CT (dashed red), empirical phase-transition curve at $50\%$ contour line (cyan), $5\%$ and $95\%$ contour lines (yellow and magenta). \label{fig:random_image_domain}}
  \end{figure}

We create DT phase diagrams as in the the previous section for the signedspikes class reconstructed by \ellone{} and spikes reconstructed by \ellp{}, see Fig.~\ref{fig:random_image_domain}. ALMT phase diagrams are omitted for brevity. As the purpose of this study is not to compare with the Gaussian sensing matrices but equi-angular fan-beam CT sampling, we do not show the theoretical phase-transition curves as in the previous section but instead in dashed red line the empirical phase-transition curves for the equi-angular fan-beam CT geometry, which was shown in cyan in Fig.~\ref{fig:dtalmt_signedspikes} and Fig.~\ref{fig:dtalmt_spikes}.

Compared to the equi-angular fan-beam case, we observe essentially no difference for the fanbeam\_rand case: The empirical phase-transition curves follow the dashed red line closely in both signedspikes with \ellone{} and spikes with \ellp{} phase diagrams. The random\_rays setup has very similar phase diagrams, but in the signedspikes case, the transition curve is slightly lower than in the equi-angular fan-beam case. In other words, on this set of image-domain sparsity test cases, randomness does not lead improved recoverability, but rather comparable or slightly reduced.

  \subsection{Gradient-domain sparsity}
  
\begin{figure}[tb]
\centering
\newcommand{\ww}{0.4\linewidth}
   \includegraphics[width=\ww]{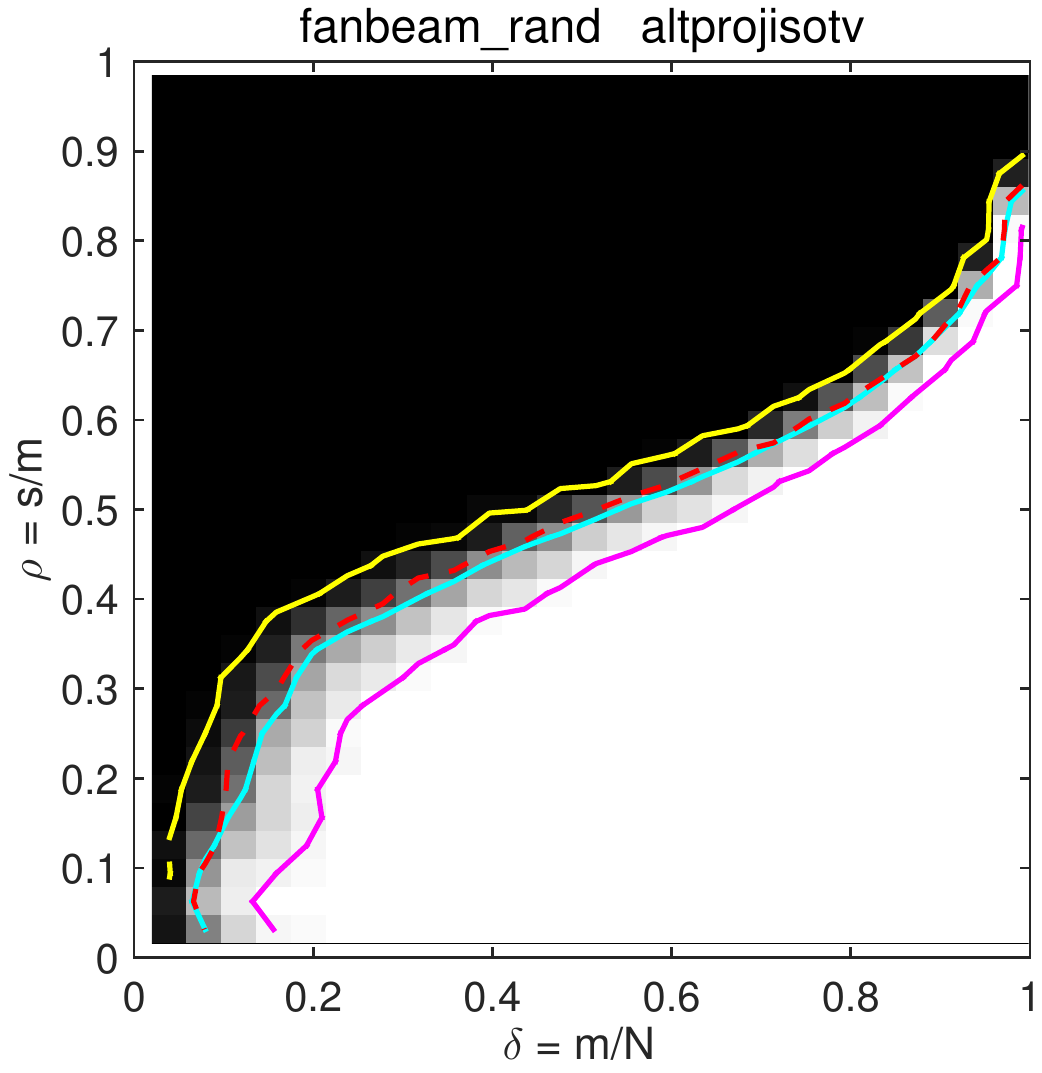}\imsep
  \includegraphics[width=\ww]{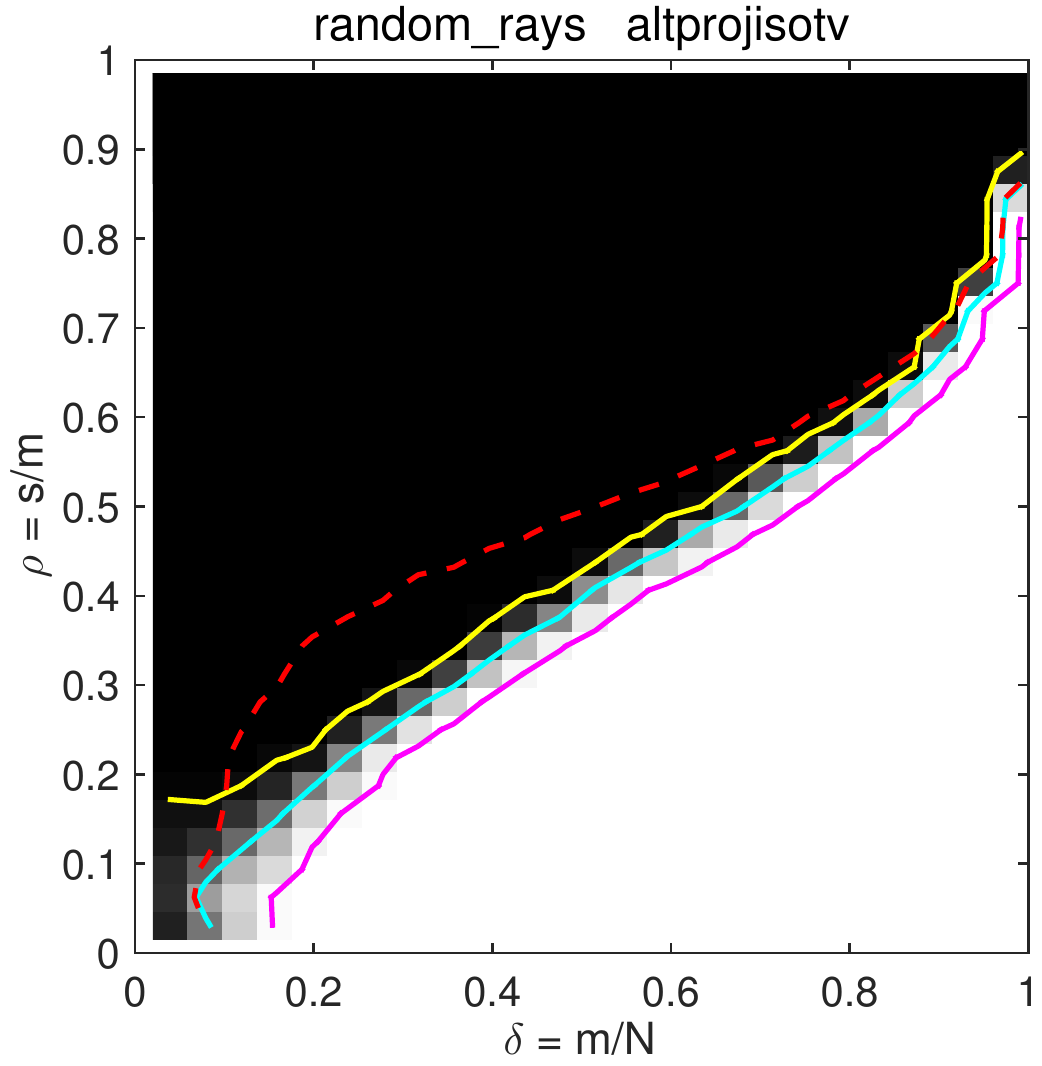}
  \caption{DT phase diagrams for the altprojisotv image class and \tv{} reconstruction. Fan-beam with random source positions (left) and random rays geometry (right). Empirical phase-transition curve for equi-angular fan-beam CT (dashed red), empirical phase-transition curve at $50\%$ contour line (cyan), $5\%$ and $95\%$ contour lines (yellow and magenta). \label{fig:random_gradient_domain}}
  \end{figure}

  For TV, we create phase diagrams for the altprojisotv class with both of the random-sampling CT setups, see Fig.~\ref{fig:random_gradient_domain}, and compare with the equi-angular fan-beam results in Fig.~\ref{fig:dtalmt_altprojisotv} indicated again by dashed red line. In both TV cases we observe \emph{worse} recoverability than equi-angular fan-beam. 
  
  The fanbeam\_rand setup has a slightly lower empirical phase-transition curve and the transition is wider than for equi-angular fan-beam, as indicated by the larger distance between the $5\%$ and $95\%$ contour lines. This means that on average slightly more projections are needed to recover the same image and further that the critical sampling level sufficient for recovery is less well-defined than for the equi-angular fan-beam case where the phase-transition is sharper. 
  
  For random\_rays the transition curve is substantially lower, meaning that on average more projections are needed for recovery of a same-sparsity image compared to equi-angular fan-beam. The largest difference is seen in the left half of the phase diagram, i.e. at fewer samples. One possible explanation of the reduced recoverability here is that with relatively few and independent rays, the probability that some pixels are not intersected by any ray is relatively large. Thus there is no information about such a pixel in the data, so the reconstructed value is solely determined by the regularizer. In contrast, in a fan-beam setup with dense projection-view sampling as in our case, all pixels will be intersected by at least one ray from each projection view.
  
\subsection{Conclusion on Study B}
By use of phase-diagram analysis we have compared two random-sampling strategies for CT with the more standard equi-angular fan-beam CT. The analysis revealed, in contrast to what might have been anticipated from the key role of randomness in CS, that random sampling does not improve recoverability in CT. On the contrary, in some cases random sampling even leads to worse recoverability, most notably for the random\_rays setup.

%%%%%%%%%%%%%%%%%%%%%%%%%%%%%%%%%%%%%%%%%%%%%%%%%%%%%%%%%%%%%%%%%%%%%
\section{Study C: Linking to realistic CT systems}
\label{sec:large-scale}

In this section, we begin the task of linking the small-scale recovery results to
realistic CT systems. 
What we are interested in is whether phase diagrams can be used to predict critical sampling levels as function of sparsity in a realistic CT system. 
The studies presented should not be regarded as complete,
and many issues for future research will be highlighted. Broadly speaking, the two
main areas of concern are test phantom and optimization algorithm.
A good test phantom presents a challenge. The small-scale phase-diagram results use phantom
ensembles generated from a probabilistic model. While the results provide
a sense of group recovery, a realization from any of the considered object
models does not look like an actual object that would be CT-scanned.

Which optimization algorithm to use is also an important question. For the small-scale studies MOSEK is a convenient choice because a highly accurate solution can be computed reliably and reasonably fast. This means that whether or not an image is recoverable can be easily verified numerically. Optimization algorithms for large-scale CT systems can not
involve more expensive operations than matrix-vector products, at present, ruling out software packages such as MOSEK in favor of first-order methods that are inherently less accurate, in particular for large-scale problems, where in practice it is often necessary to truncate iteration early. As we will show, having less accurate solutions makes it more difficult to decide whether an image is recoverable.

As large-scale studies are necessarily sparse, we cannot
provide comprehensive empirical evidence of sufficient sampling but only a preliminary indication
of how well phase-diagram analysis can predict sufficient sampling for SR for realistic CT systems. As we will show, even this is a complex task for example due to complicated image structure and algorithmic issues, and we will point out several future directions to pursue.

Sec.~\ref{sec:ls-phantom} presents two phantoms generated for
the present study to have different levels of realism with respect to an actual CT scanned object. 
Sec.~\ref{sec:ls-alg} presents the first-order optimization algorithm we use for the large-scale
recovery studies, while Sec.~\ref{sec:ls-issues} illustrates some of the algorithmic and numerical challenges we face. 
Sec. \ref{sec:ls-results} shows recovery results for
the two phantoms as a function of number of CT projections and compare with critical sampling levels predicted from small-scale phase diagrams.

%%%%%%%%%%%%%%%%%%%%%%%%%%%%%%%%%%%%%%%%%%%%%%%%%%%%%%%%%%%%%%%%%%%%%
\subsection{Walnut test phantoms}
\label{sec:ls-phantom}

\begin{figure}[tb]
\centering
\newcommand{\ww}{0.65\linewidth}
\includegraphics[width=\ww]{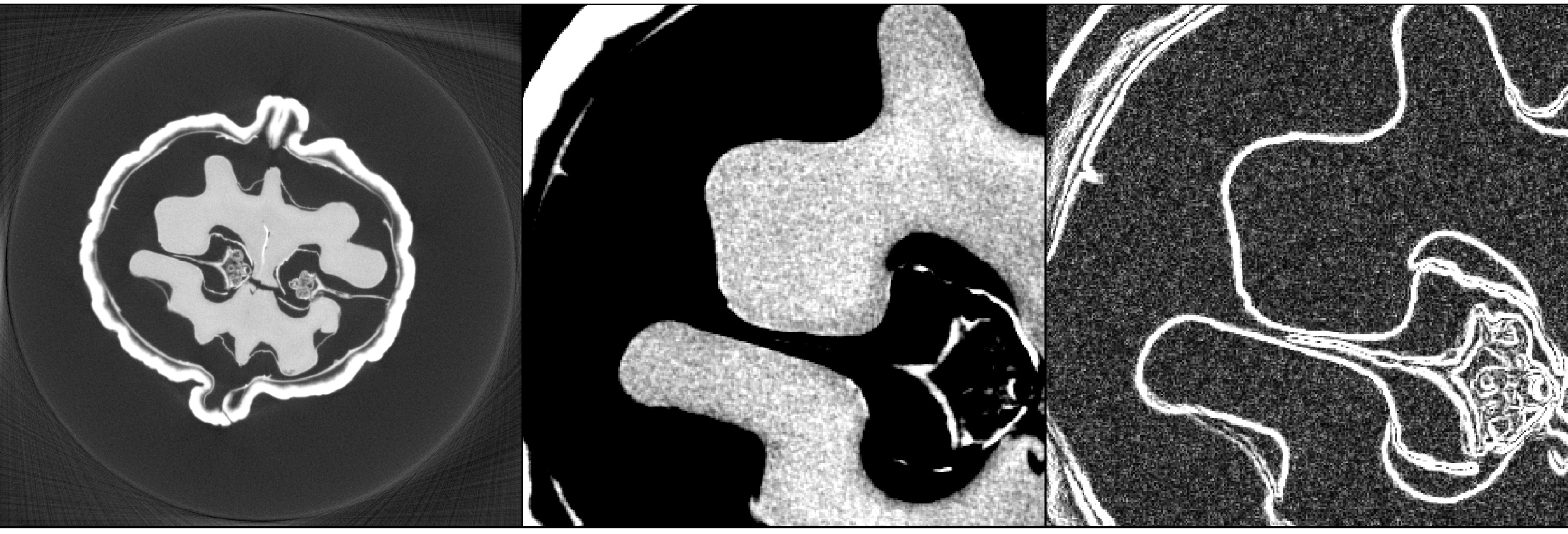}\\
\includegraphics[width=\ww]{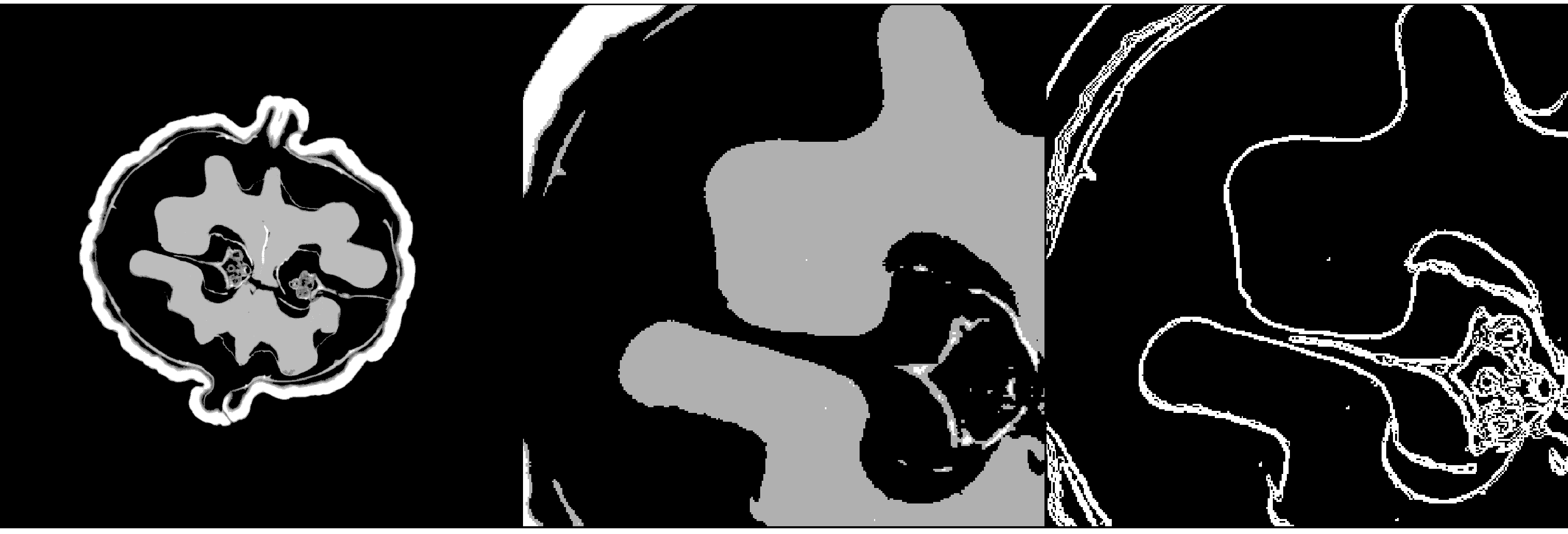}\\
\includegraphics[width=\ww]{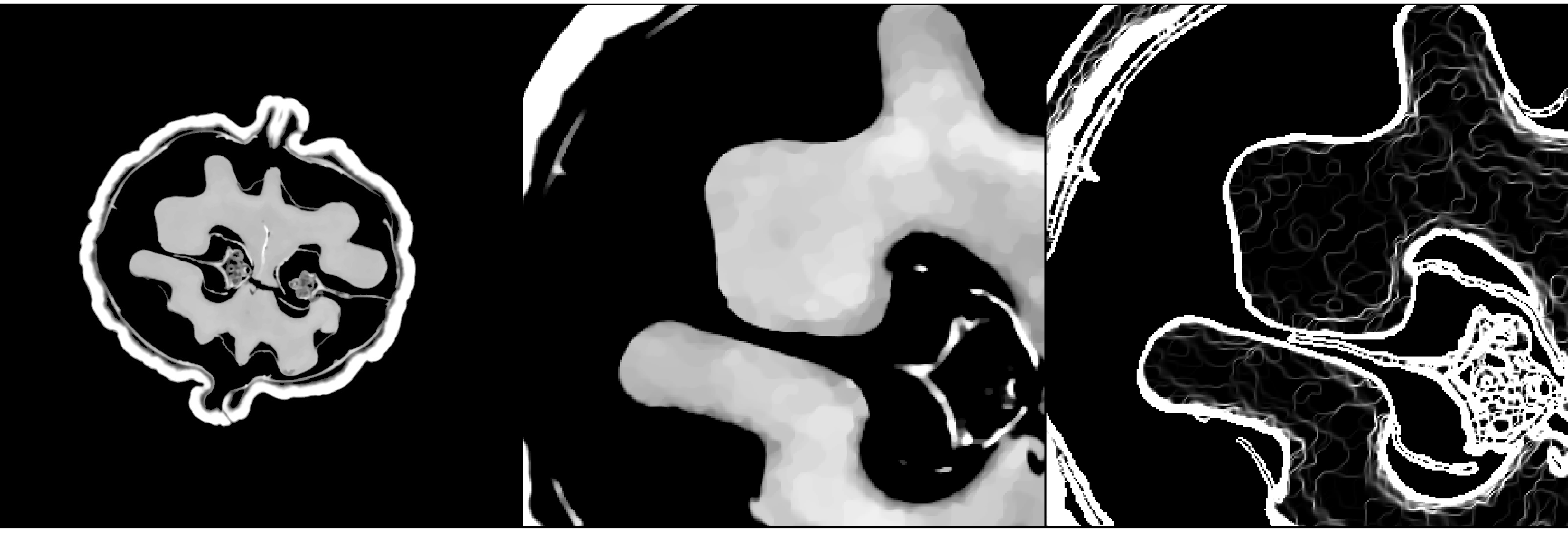}
\caption{(Top row) tomographic slice of a walnut, (middle row) structure phantom derived
from the walnut slice image, and (bottom row) texture phantom also derived from this image.  The left column shows the whole image in the gray scale window of
[0,0.5] cm$^{-1}$, except for the original walnut image where it is [-0.1,0.5] cm$^{-1}$.
The middle column shows a blown-up region of interest in the narrower gray scale window
[0.3,0.4] cm$^{-1}$ in order to see the texture on the walnut meat. The right column illustrates the gradient-magnitude image in the gray scale window [0,0.01] cm$^{-1}$,
except for the original walnut image where it is [0,0.05] cm$^{-1}$.
\label{fig:walnutScan}}
\end{figure}

In the present large-scale study, there are two links that need to be established to
relate the small-scale phase-diagram analysis to realistic CT: the system size
needs to be extrapolated up, in this case to $N_\text{side}=1024$; and the results from
the various probabilistic phantom models need to extend to realistic structure as seen
in actual CT-scan objects. We address both by designing two large-scale test phantoms with increasing realism from an actual CT scan of a walnut. The idea of scanning a walnut comes from \cite{SiltanenTV:2014}.

In choosing a test phantom for image recovery studies, we aim for an image
with gradient-domain sparsity to illustrate the effectiveness of TV in reducing the necessary
number of samples for accurate image recovery. Yet, the phantom should also have features
somewhat representative of what would be encountered in CT applications.
Typical computer phantoms for CT image
reconstruction testing, composed of simple geometric shapes of uniform gray levels,
are unrealistically sparse in the gradient domain. Such phantoms would be helpful in 
extrapolation of small-scale phase-diagram analysis, but  do not have much bearing in actual CT applications. 

The basis of the test phantoms we generate is a cone-beam CT scan data set of a walnut. The data consists of 1600 equiangular $1024^2$-pixel projections acquired on a Zeiss Xradia 410 Versa micro-CT scanner operated at a 40 kV source voltage, 5 s exposure per projection, and $10.51$ cm source-to-center and $4.51$ cm center-to-detector distances. The central slice is reconstructed onto a $1024^2$-pixel image (pixel size $46.0773 \cdot 10^{-6}$ m) from the corresponding rows of data using 500 iterations of a SIRT-type algorithm.

The first and simplest phantom, the \emph{structure} phantom, is derived from this image
by equalizing the image gray value histogram to 7 discrete gray levels including the background value of 0.
The second and more complex phantom, the \emph{texture} phantom, is derived from the walnut image by performing TV-denoising on the original walnut image after thresholding small background pixel values to zero. 
The two versions of the walnut phantoms including blow-ups and gradient-domain images are shown in Fig.~\ref{fig:walnutScan} and gradient-domain sparsity values are given in Table~\ref{tab:walnutdata}.
The studies are idealized in that there is no data
inconsistency; in actual CT the projection data $b$ will in general not be in the range
of the projection operator $A$, and there is in this case no solution to the linear
system $Ax=b$.

%%%%%%%%%%%%%%%%%%%%%%%%%%%%%%%%%%%%%%%%%%%%%%%%%%%%%%%%%%%%%%%%%%%%%
\subsection{Large-scale first-order optimization algorithm}
\label{sec:ls-alg}
\begin{algorithm}[tb]
\hrulefill
\begin{algorithmic}[1]
\State INPUT: data $b$
\State INPUT: tuning parameter $\lambda$
\State $\nu = \| A \|_2 / \| S \|_2$
\State $L \gets \| (A, \nu S) \|_2$
\State $ \tau \gets 1/L; \; \sigma \gets 1/L; \; \theta \gets 1; \; k \gets 0$
\State initialize $x_0$,  $y_0$, and $z_0$ to zero vectors
\State $\bar{x}_0 \gets x_0$
\Repeat
\State $y_{k+1} \gets  y_k+\sigma( A \bar{x}_k -b)$ \label{dualDataUpdate}
\State $ z_k^\prime = z_k + \sigma \nu S \bar{x}_k$
\State $z_{k+1} \gets z_k^\prime 
( (\lambda / \nu) / \max(\lambda /\nu, | z_k^\prime |))$ \label{dualGradUpdate}
\State $x_{k+1} \gets  x_k - \tau (\sm{A}^T y_{k+1} +
\nu S^T z_{k+1}) $ \label{primalupdate}
\State $\bar{x}_{k+1} \gets x_{k+1} + \theta(x_{k+1} - x_k)$
\State $k \gets k+1$
\Until{$k \ge K$}
\State OUTPUT: $x_K$
\end{algorithmic}
\hrulefill
\caption{Pseudo-code for $K$ steps of the CP algorithm instance
for solving Eq. (\ref{p1tv}). When $S=I$ and $S=D$ this algorithm applies to
\ellone{} and \tv{}, respectively. The variables $y_k$ and $z_k$ are dual to the sinogram
and 
image, respectively. For gradient-domain sparsity (TV) $z_k$ has
the dimension of the image gradient, and for image-domain sparsity (\ellone{})
$z_k$ has the dimension of the image itself.}
\label{alg:p1tv}
\end{algorithm}

We consider large-scale solvers for problems \ellone{} and \tv{}. There has been much
recent research on first-order algorithms \cite{beck2009fast,chambolle2011first},
motivated by exactly the type
of problem we face here. We require a solver that can handle the non-smoothness
of \ellone{} and \tv{}, and which can be applied to large-scale systems such as CT,
where the images can contain 10$^6$ pixels in 2D or 10$^9$ voxels in 3D and data sets of similar size.
The CT system specifically presents another challenge in that the system
matrix representing standard X-ray projection has poor conditioning \cite{jakob:2013}.
An additional difficulty in solving \ellone{} and \tv{}, compared to the form \eqref{eq:ineqreconprob}, is in satisfying the equality constraint;
achieving this constraint to numerical precision with present computational and algorithmic
technology is not possible as far as we know. We present, here, our adaptation
of the Chambolle-Pock (CP) primal-dual algorithm, which we have found to be effective
for the CT system \cite{sidky:CP:2012,sidky2013first,SidkyJTEHM:2014}.

The algorithm used is essentially the same as the one developed in Ref. \cite{SidkyJTEHM:2014}.
The CP algorithm instance is designed to solve the following optimization problem
\begin{equation}
\label{p1tv}
\argmin_x ~ \frac{\lambda}{\nu} \sum_j \|  \nu S_j x \|_2 \quad \text{subject to} \quad Ax = b,
\end{equation}
where Eq. (\ref{p1tv}) becomes \ellone{} and \tv{} when the sparsifying operator is $S_j=I_j$ and
$S_j=D_j$, respectively; $I_j$ is an image where the $j$th pixel is one and all other pixels are zero; $\nu$ is a constant which balances the operator norms
\begin{equation}
\notag
\nu = \| A \|_2 / \| S \|_2,
\end{equation}
where $S$ is a matrix of $S_j$ for all $j$;
and the parameter $\lambda$, which does not affect the solution of Eq. (\ref{p1tv}),
is used to improve numerical convergence. The parameter $\lambda$ is tuned empirically.
The corresponding algorithm for solving Eq. (\ref{p1tv})
is shown in pseudo-code form in Alg.~\ref{alg:p1tv}.

Considering that the phantom-recovery studies we want to use Alg.~\ref{alg:p1tv} for involve multiple runs over different system
matrices $A$ corresponding to CT sampling with different numbers of projections, we
found it most practical to obtain results for fixed iteration number $K$ and tuning
parameter $\lambda$.  The computational time for performing the expensive operations
$Ax$ and $A^T x$ makes consideration of a prescribed stopping criterion difficult.
For the $N_\text{side}=1024$ system of interest these time-limiting operations take
~1 second for our GPU-accelerated projection codes.
A fixed stopping criterion entails variable numbers of iterations, and we
have observed that for Alg. \ref{alg:p1tv} the number of iterations can vary from
1,000 to over 100,000 iterations for a convergence criterion of interest. In terms of
computational time, this range translates to ~20 minutes to well over a day.
As a result, a study may not be completed in a reasonable amount of time;
thus, we fix $K$ and $\lambda$ for our phantom-recovery study. 

Because 
large-scale first-order optimization algorithms are seeing many new developments at present, it is
likely that there either exists or will be a better alternative to Alg.~\ref{alg:p1tv}.
In fact, we invite the interested reader to find such an alternative, which
can have an important impact on CT imaging!
For example, as will be seen shortly, Alg.~\ref{alg:p1tv} 
has limited success for phantom recovery studies for \ellone{} on systems of realistic size.

\subsection{Algorithm issues} \label{sec:ls-issues}

We demonstrate first some of the challenges in carrying out large-scale recovery studies by applying Alg.~\ref{alg:p1tv} to a medium-scale problem 
using a $N_\text{side}=128$ version of the structure walnut phantom.
The phantom has a gradient-domain sparsity of 1826 and a pixel sparsity of 2664 out of a total of 11620 pixels.
We do a recovery study by studying the root-mean-square (RMSE) reconstruction error as function of the number of projections.
We discuss in detail specific issues of the sampling recovery study
for the purpose of understanding the large-scale results. 

\begin{figure}[tb]
\centering
 \newcommand{\ww}{0.49\linewidth}
 \includegraphics[width=\ww]{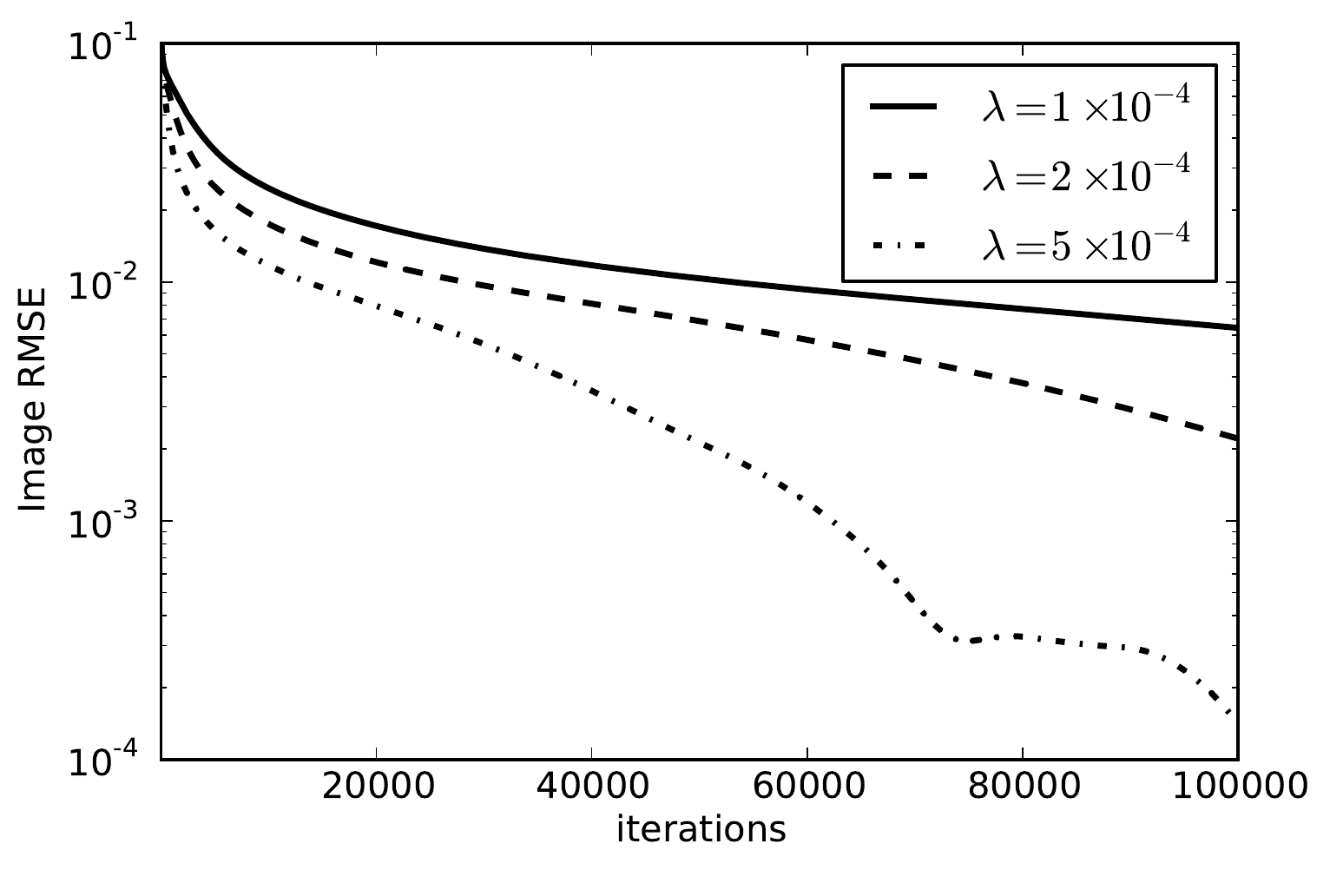}
 \includegraphics[width=\ww]{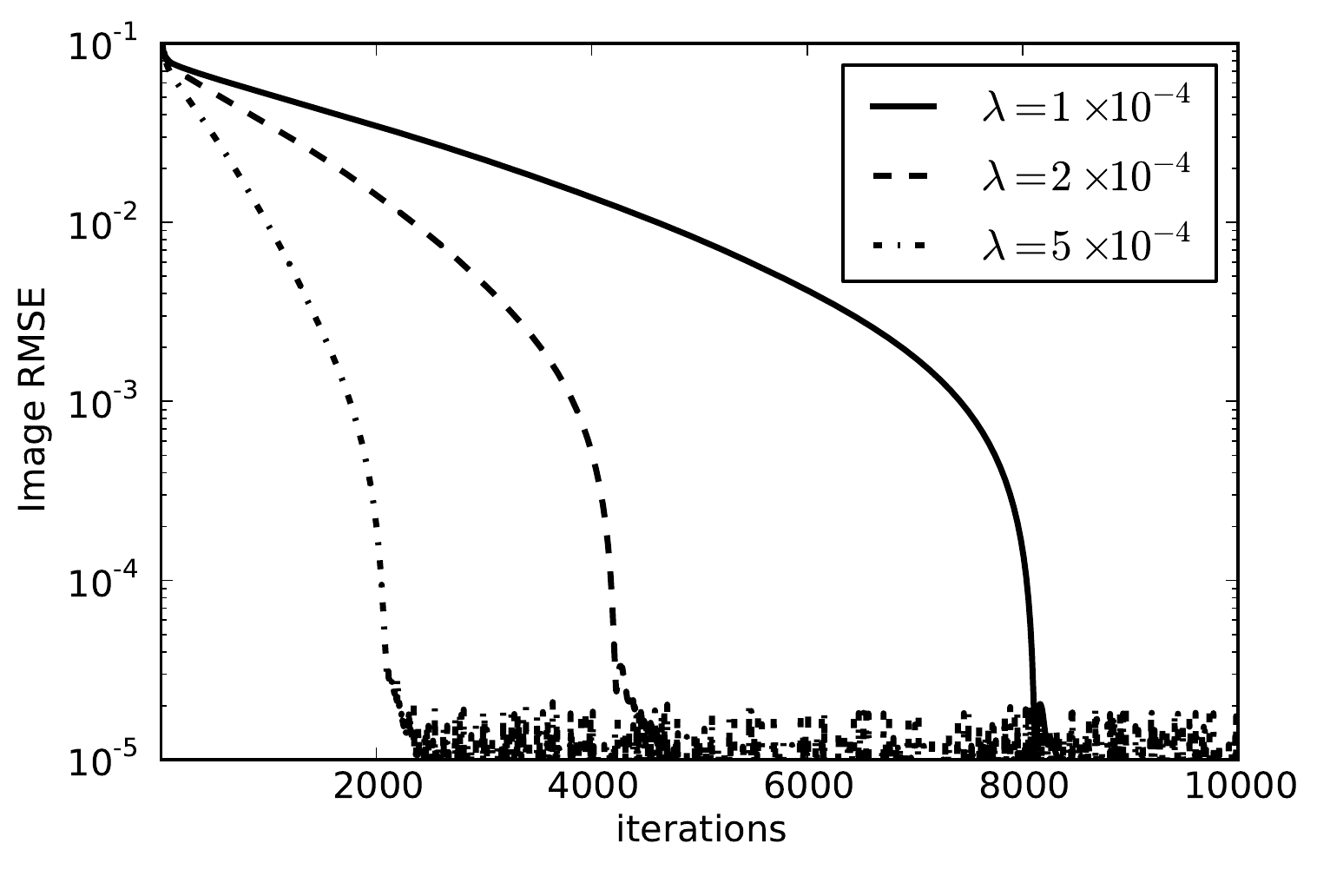}
  \caption{Image RMSE curves resulting from Alg. \ref{alg:p1tv} run with different values of $\lambda$ for $N_\text{v}=21$ and data generated from the $N_\text{side}=128$ version of the structure walnut phantom. Results for \ellone{} and \tv{} are shown on the left and right, respectively. \label{fig:lambdaconv}}
\end{figure}

\paragraph*{The tuning parameter $\lambda$ and convergence}
The tuning parameter $\lambda$ does not affect the solution of \ellone{} or \tv{}, but it can
have a large impact on convergence. To illustrate this we show results of single runs for
$N_\text{side} = 128$ and $N_\text{v}=21$ for both \ellone{} and \tv{} in Fig.~\ref{fig:lambdaconv}. 
The value $N_\text{v}=21$ is chosen because it is the smallest number of views for which accurate recovery
is obtained for both \ellone{} and \tv{}.
Note that we are showing results for $K=100,000$ iterations for \ellone{}, while only $K=10,000$ for \tv{}.
\begin{figure}[tb]
\centering
 \newcommand{\ww}{0.49\linewidth}
 \includegraphics[width=\ww]{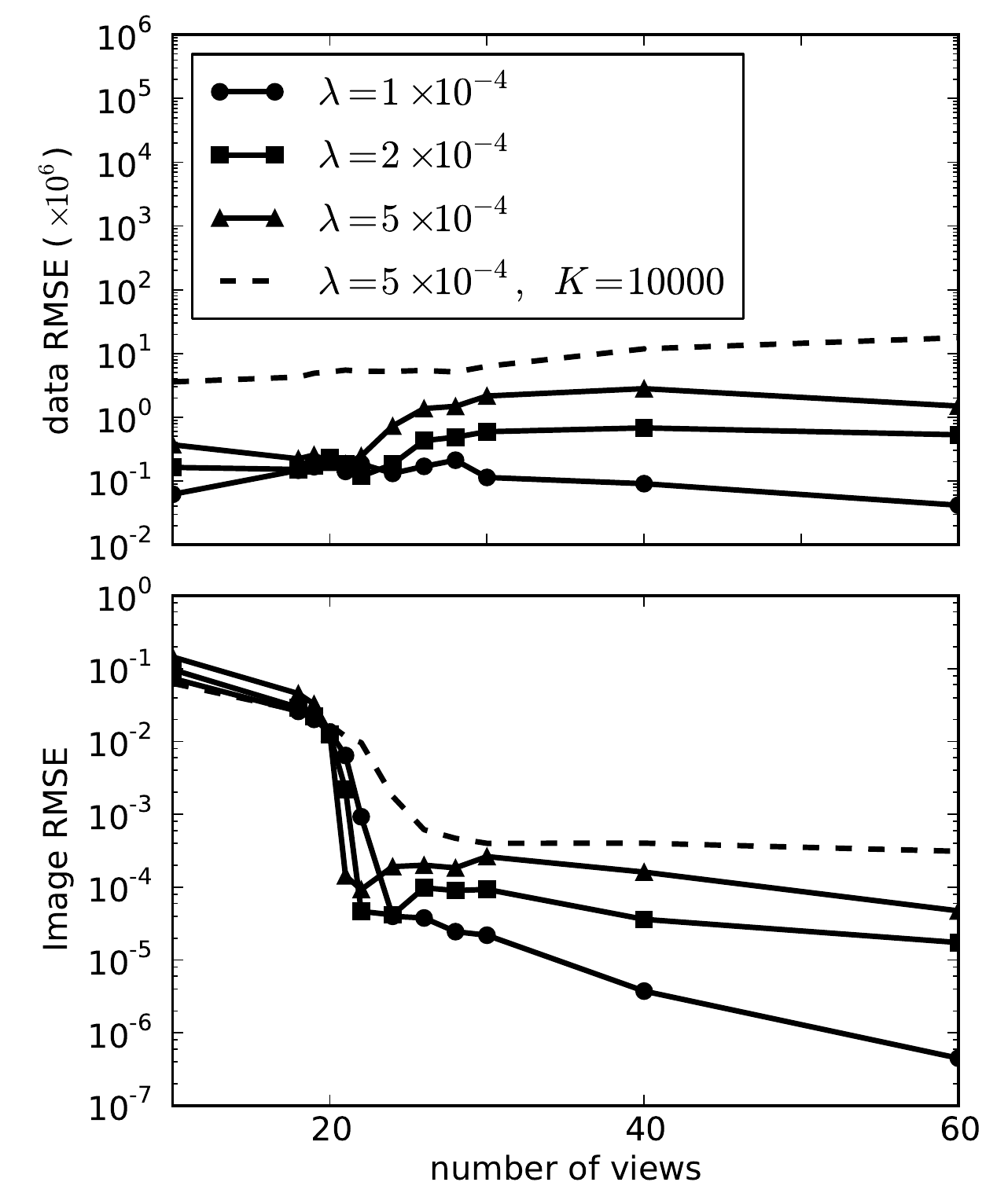}
 \includegraphics[width=\ww]{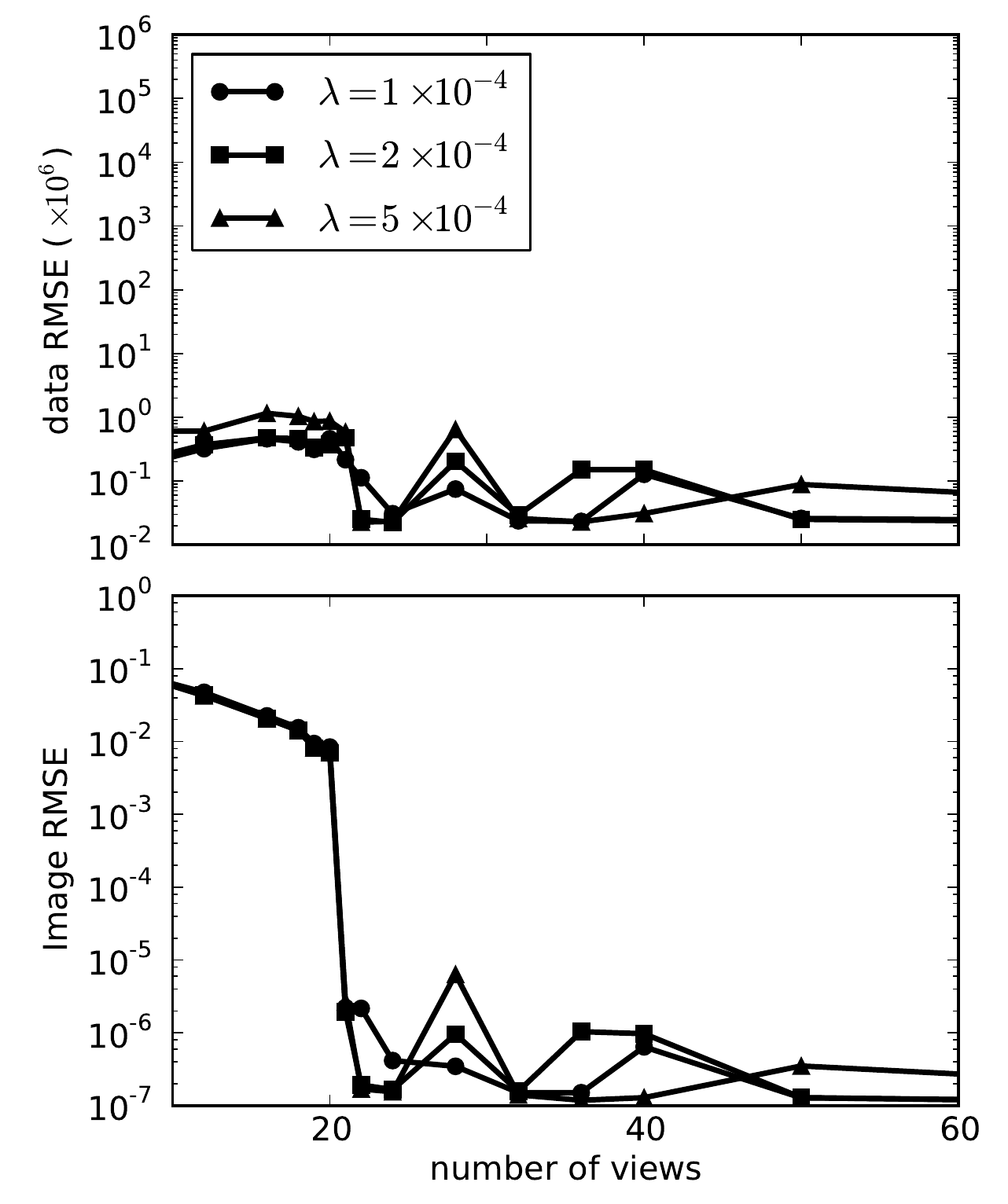}
 \caption{Image and data RMSE plots for the $N_\text{side}=128$ version of the structure walnut phantom using Alg. \ref{alg:p1tv} with different values of $\lambda$. The results for \ellone{} (left) are obtained for $K=10^5$ iterations except for the indicated curve for $K=10^4$. The results for \tv{} (right) are obtained for $K=10^4$ iterations.\label{fig:lambda}}
\end{figure}
It is clear that convergence rates change significantly with $\lambda$, and consequently
recovery curves will be affected by $\lambda$. While $\lambda$ is specific to Alg.~\ref{alg:p1tv}, optimization algorithms generally entail parameters with large effect on convergence rate.

The impact on recovery curves is seen in Fig.~\ref{fig:lambda}, where we
we compare recovery curves obtained at different $\lambda$ for \ellone{} ($K=100,000$ iterations) and \tv{} ($K=10,000$ iterations). While overall the recovery curves are similar some differences appear in particular near the jump in error for \ellone{}. This can complicate the accurate estimation of the jump location. Overall, in this case, the lowest
image RMSE is obtained for $\lambda=5\times10^{-4}$. For the large-scale system $N_\text{side}=1024$ we have found the value of $\lambda=1\times10^{-4}$ to be useful
for \ellone{} and \tv{}, and for different values of $N_\text{v}$ and $N_\text{side}$. One could envision
a strategy where Alg.~\ref{alg:p1tv} is run with a small set of $\lambda$ values and the lowest
image RMSE at iteration $K$ is taken for the recovery plot. In the large-scale results presented shortly,
we found this to be unnecessary, and $\lambda$ is simply fixed at $1 \times 10^{-4}$.

\paragraph*{Recovery plots and difficulty with \ellone{}}
The phantom recovery plots for \ellone{} and \tv{} in Fig.~\ref{fig:lambda} both show the distinct jump in RMSE at a certain number of projections, at which the image is recovered. We recognize this from the small-scale $N_\text{side}=64$ studies in \cite{Joergensen_eqconpap_v2_arxiv:2014}. 
The price of using fixed $K$, however, is that convergence results across projection numbers
are not uniform as the data discrepancy varies with view number.

Furthermore, the recovery curve can be severely affected by poor convergence. If instead of $K=100,000$ only take $K=10,000$ as in the \tv{} case, the remaining recovery curve in Fig.~\ref{fig:lambda} is obtained.  The previously abrupt change in error is considerably smoothed and shifted to a different number of views.

The issue of convergence, here, is ubiquitous in iterative image reconstruction
for CT and it can be traced to the use of matched projection, $A$, and back-projection,
$A^T$, where it is well-known in the CT community that matched projector/back-projector
pairs can lead to Moire artifacts that decay extremely slowly \cite{de2004distance}. As a result,
many iterative algorithms in CT employ a different back-projection matrix $B \ne A^T$,
\cite{zeng2000unmatched}. For our purpose we must use the matched pair, in order to solve a well-defined optimization problem.
For the larger system, sufficient iteration for \ellone{} lies out of reach with Alg.~\ref{alg:p1tv} and we
focus only on phantom recovery for TV.

\subsection{Large-scale recovery results}
\label{sec:ls-results}

\paragraph*{Predicting sufficient sampling from phase diagrams}
\begin{figure}[tb]
\centering
\includegraphics[width=0.4\linewidth]{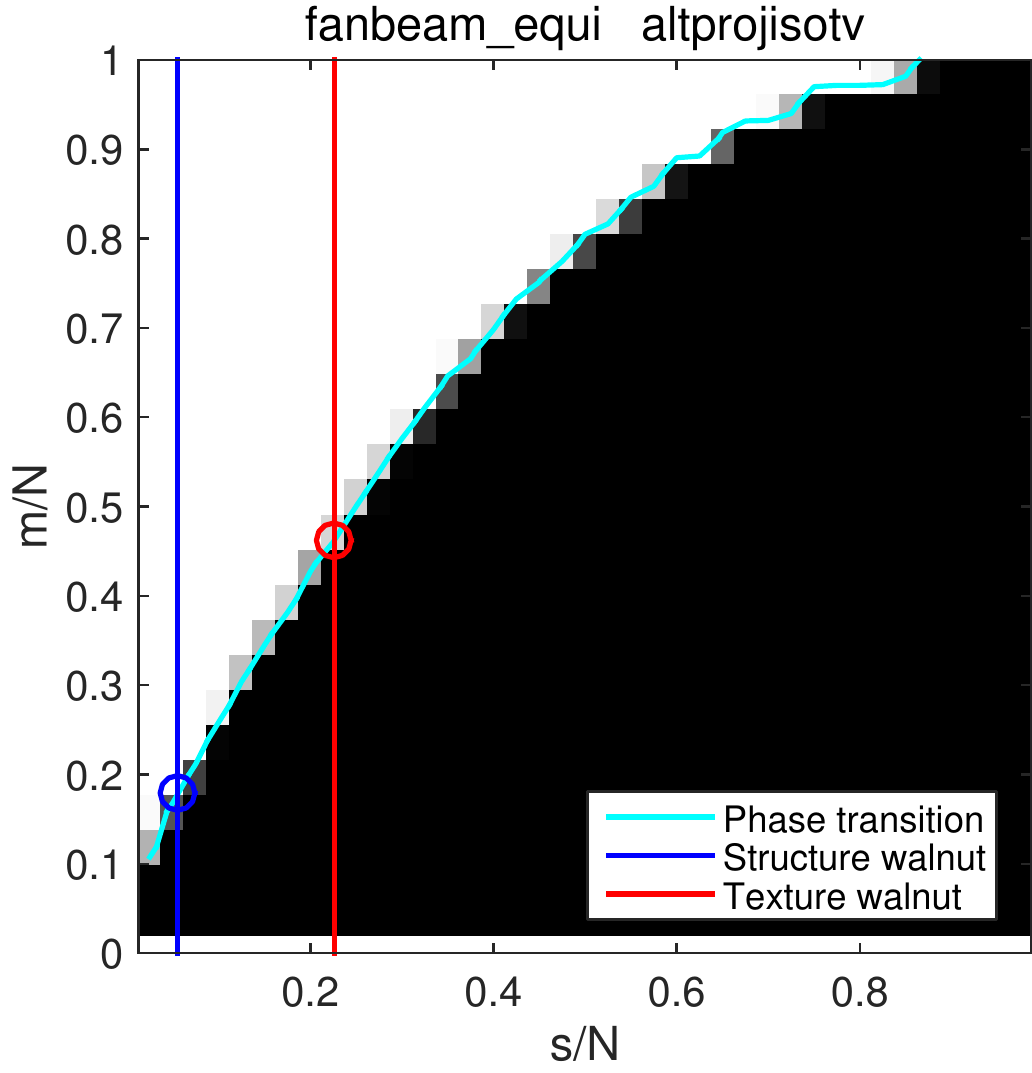}\imsep
    \includegraphics[width=0.405\linewidth]{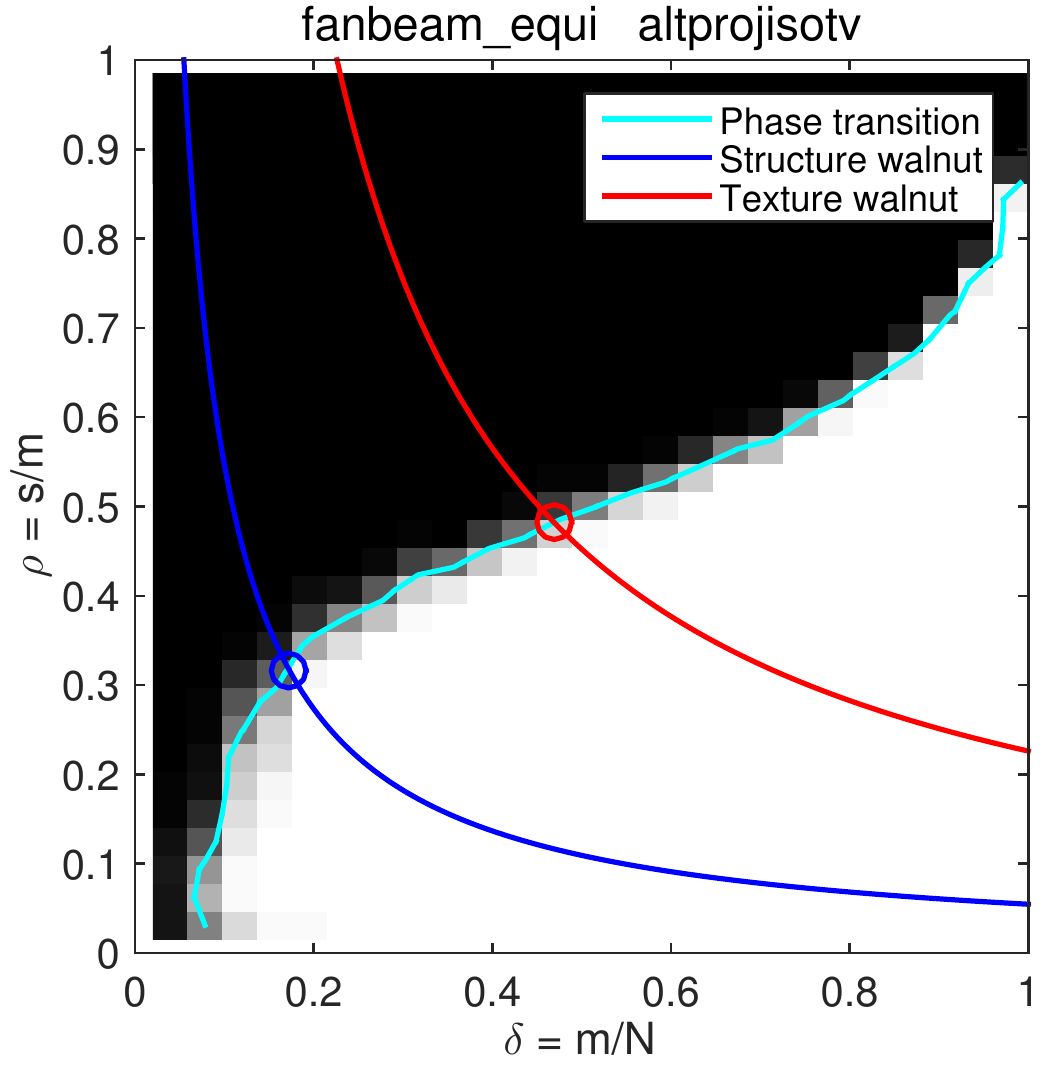}
      \caption{Prediction of critical sampling for TV and walnut phantoms by ALMT (left) and DT (right) phase diagrams. \label{fig:prediction}}
\end{figure}
We will use the phase diagrams from Study A to predict critical sampling levels for large-scale \tv{} reconstruction. We found in \cite{Joergensen_eqconpap_v2_arxiv:2014} that the ALMT phase diagram of a given image class remains unchanged 
at image resolutions $N_\text{side} =32$, $64$ and $128$, i.e., is independent of resolution. We assume this holds also for the DT phase diagram and we use the DT and ALMT phase diagrams from Fig.~\ref{fig:dtalmt_altprojisotv} (which are for $N_\text{side}=64$) to predict critical sampling levels for the two walnut images at $N_\text{side} = 1024$.

We illustrate in Fig.~\ref{fig:prediction} how to determine critical sampling levels given a sparsity level. 
The number of pixels inside the disk is $823592$ and the gradient sparsity levels of the structure and texture walnut images are given in Table~\ref{tab:walnutdata}. In the ALMT phase diagram we can trace vertical lines at each $s/N$ value and find the intersections (indicated by circles) with the empirical phase-transition curve, which gives the predicted critical $m/N$ values. By multiplication of $N$ and division by the number of rays in a single projection, i.e. $2048$, we get the critical number of projections, see Table~\ref{tab:walnutdata}.

To do the same in the DT phase diagram we combine $\delta = m/N$ and $\rho = s/m$ into $\rho = s/(\delta N)$, i.e., a fixed sparsity $s$ traces out a hyperbola on $\delta \in (0,1)$. For the hyperbola of each walnut image we find the intersection point $(m/N, s/m)$ with the empirical phase-transition curve. Up to the accuracy of reading off the figure, the two components lead to identical critical values of $m$, from which we find the critical number of projections for each walnut image, see Table~\ref{tab:walnutdata}.

We note that the larger number of gradient non-zeroes in the texture walnut image leads to prediction of a higher critical sampling level. Similar plots for image-domain sparsity could be constructed based on Fig.~\ref{fig:dtalmt_signedspikes} and Fig.~\ref{fig:dtalmt_spikes} and the fixed-sparsity curves would then reflect that the walnut images have more non-zeroes in the pixel domain than in the gradient domain, yielding higher predicted critical sampling levels for \ellone{}/\ellp{} than for \tv{}.

\paragraph*{Recovery of the large-scale walnut phantoms}
\begin{figure}[tb]
\centering
 \newcommand{\ww}{0.49\linewidth}
 \includegraphics[width=\ww]{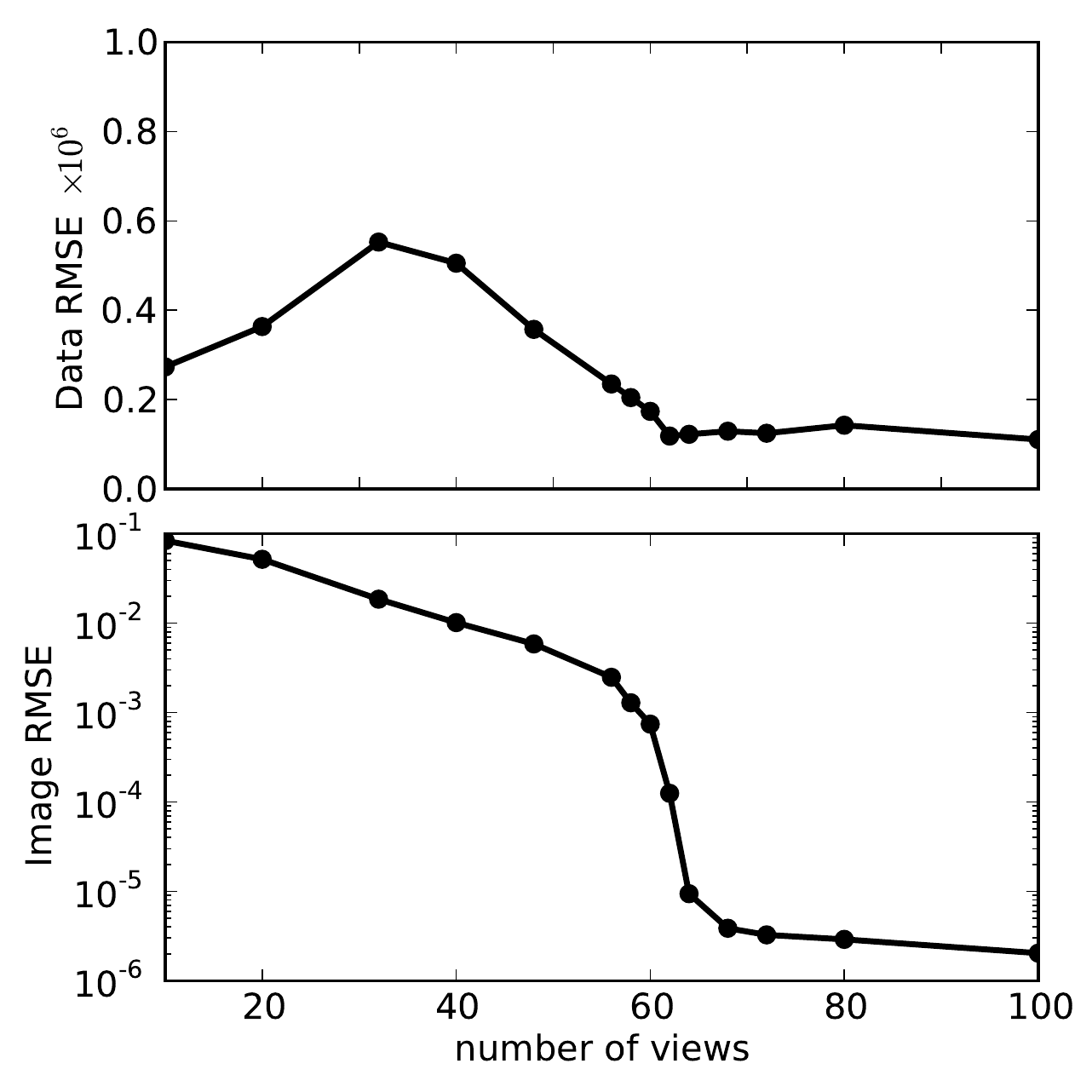}
 \includegraphics[width=\ww]{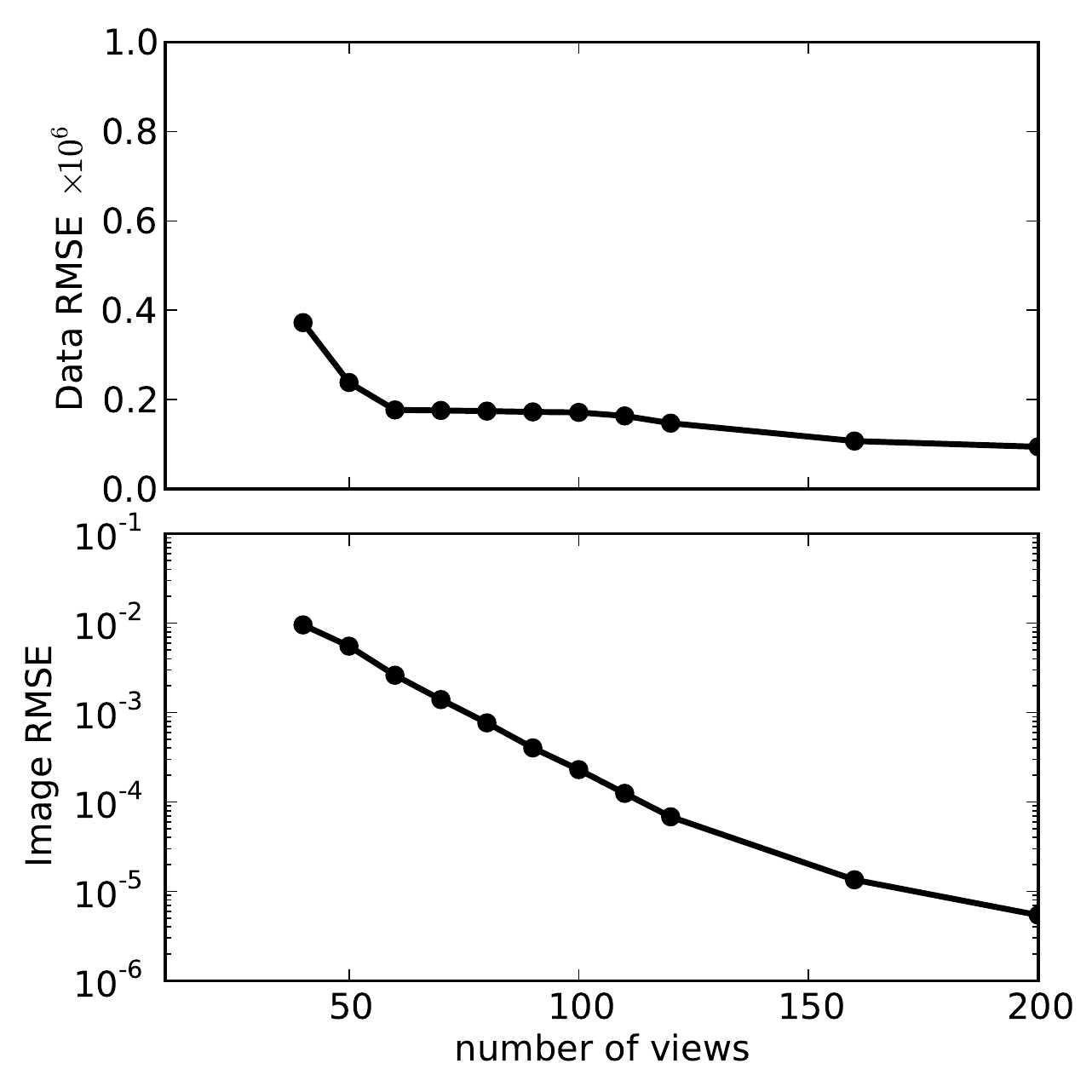}
 \hspace*{\ww}
 \caption{
Image and data RMSE plots for the $N_\text{side}=1024$ structure (left) and texture (right) walnut phantom using Alg. \ref{alg:p1tv} with $\lambda=10^{-4}$. The results are obtained at $K=10^4$ iterations.
\label{fig:walnuterrorcurves}}
\end{figure}
\begin{table}[htb]
\centering
\begin{tabular}{l|c|c|c|c}
Walnut image & Gradient sparsity & Recovered at & DT prediction & ALMT prediction \\
\hline
Structure & $\phantom{1}45,074$ & $68$ & $\phantom{1}69.3$ & $\phantom{1}71.7$ \\
Texture & $186,306$ & ? & $188.7$ & $185.8$
\end{tabular}
\caption{Walnut test images with gradient-domain sparsity levels, number of projections at which recovery is observed, and DT and ALMT phase-diagram predictions of critical sampling levels. A reference point of full sampling is $\nv \geq 403$ projections, where the system matrix has more rows than columns.\label{tab:walnutdata}}
\end{table}
We employ Alg.~\ref{alg:p1tv} to solve TV on the large-scale $N_\text{side}=1024$ CT system for the structure and texture walnut phantoms.
The resulting recovery plots are shown in 
Fig.~\ref{fig:walnuterrorcurves}.
For the structure walnut we observe an abrupt change in image
RMSE with $\nv= 68$ yielding accurate recovery, as decided by the first point where there is essentially no further decrease in RMSE. The predicted critical sampling levels from DT and ALMT phase-diagram
analysis are only slightly higher at $\nv = 69.3$ and $\nv = 71.7$, respectively, cf. Table~\ref{tab:walnutdata}. This result is rather remarkable
in that the extrapolation is extended quite far from the size of the original phase-diagram
analysis. Also, the structure phantom is clearly different from any expected
realization of any of the studied probabilistic phantoms models.

The recovery curve for the texture phantom, on the other hand, does not exhibit an abrupt change in reconstruction error, rather a gradual improvement all the way up to ~200 projections.
We therefore can not point to a specific critical sampling level.

\paragraph*{Reconstructed images for the structure and texture phantoms}
\begin{figure}[htb]
\centering
\newcommand{\ww}{0.65\linewidth}
 \includegraphics[width=\ww]{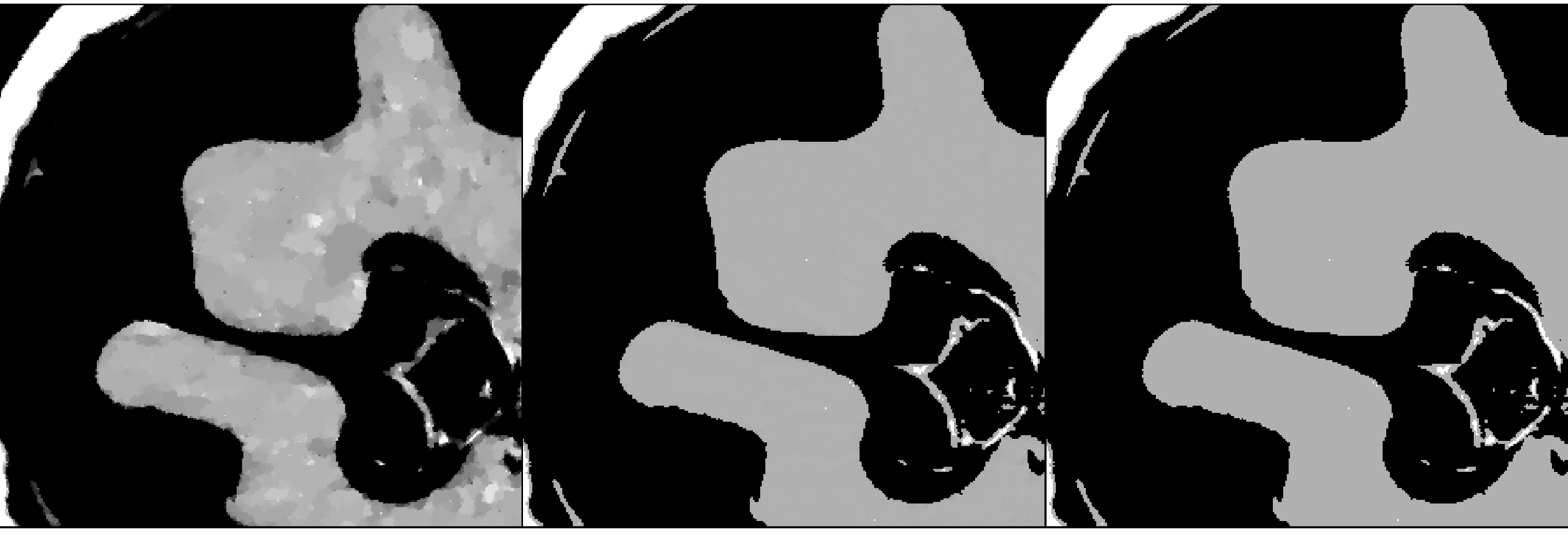}\\
 \includegraphics[width=\ww]{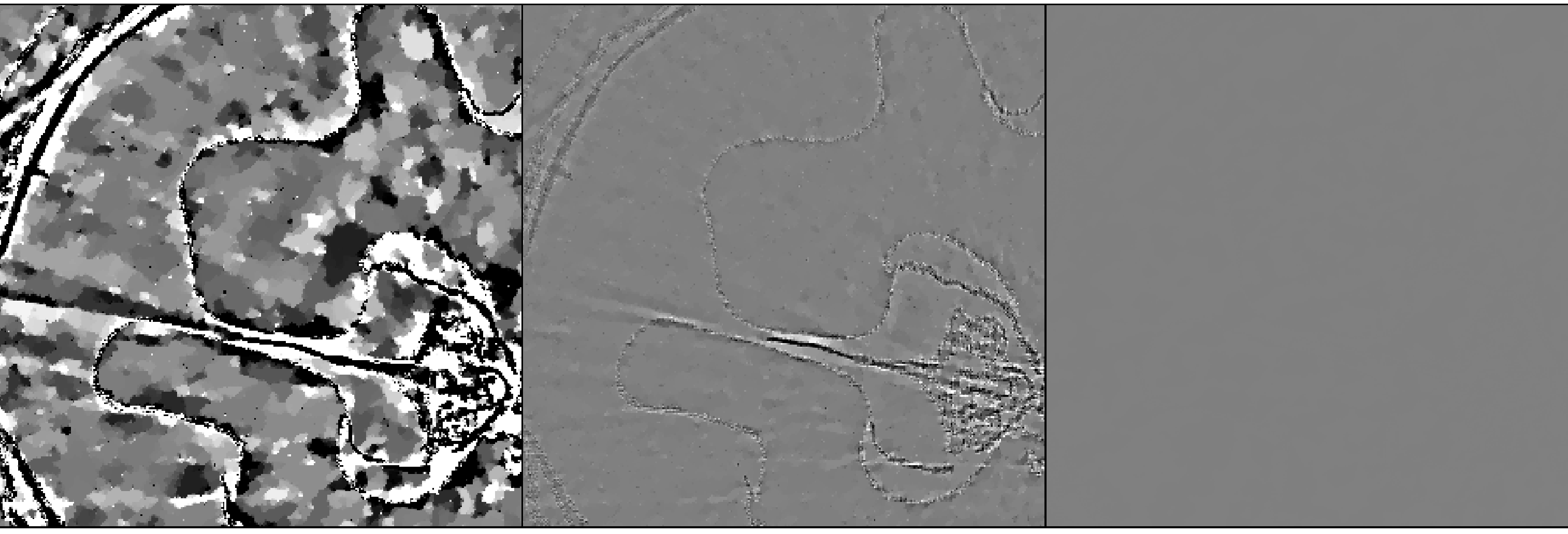}\\
 \includegraphics[width=\ww]{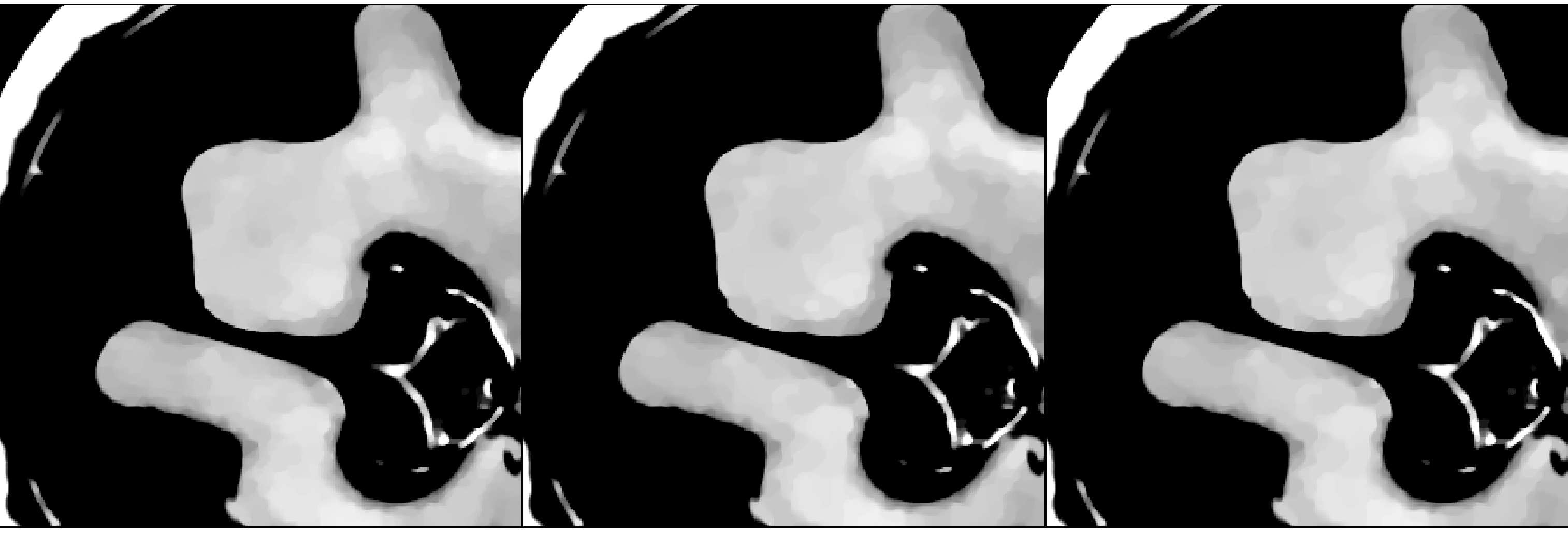}\\
 \includegraphics[width=\ww]{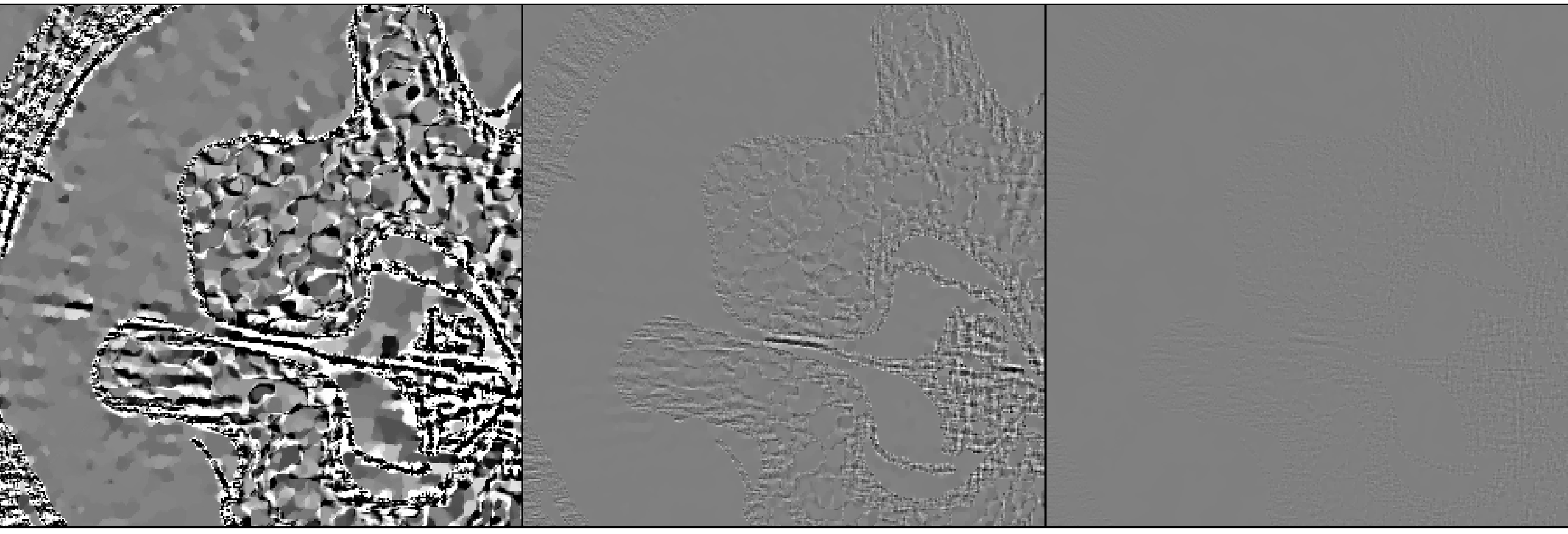}
 \caption{First row: reconstructed images from data generated by the structure
walnut with 40 (left), 60 (middle), and 68 (right) projection views (gray scale window
[0.3,0.4] cm$^{-1}$). Second row: same as first row except the structure walnut image is subtracted from the reconstructed images (gray scale window [-0.01,0.01] cm$^{-1}$).
Third row: reconstructed images from data generated by the texture
walnut with 80 (left), 120 (middle), and 160 (right) projection views (gray scale window
[0.3,0.4] cm$^{-1}$).  Fourth row: same as third row except the texture walnut image is subtracted from the reconstructed images (gray scale window [-0.001,0.001] cm$^{-1}$).
\label{fig:recons}}
\end{figure}
It is illuminating to inspect some of the reconstructed images in Fig. \ref{fig:recons}, which
correspond to the plots in Fig. \ref{fig:walnuterrorcurves}. The second and third reconstructions
for the structure phantom straddle the sharp transition in the corresponding image RMSE curve,
and it can be seen clearly in the difference image that the result for $N_v=60$ is not recovered,
while that of $N_v=68$ is much closer to the test phantom.  We point out, however, that the
difference images are displayed in a narrow 4\% gray scale window and visually the $N_v=60$ image
appears the same as the structure phantom. That the discrepancies between reconstruction and phantom
are so small emphasizes the challenge for the large-scale optimization algorithms; for actual application
where images are presented for visual inspection such
accurate solution to Eq. (\ref{eq:ineqreconprob}) would not be necessary. The results
for the texture phantom are also quite interesting in that we see the reconstructed image is visually
accurate for as few views as $N_v=80$. That there is no sharp recovery transition for the texture phantom
is likely due to the fact that the object variations occur on two scales: the jumps of the structure
borders, and the splotches of the walnut meat texture. It also can not be ruled out that a sharper recovery transition will occur if the accuracy of the computed solutions is improved even further.

\subsection{Conclusion on Study C}
In this study we have taken first steps toward phase-diagram analysis for prediction of critical sampling levels for realistic CT systems. Both test phantom design and accurate large-scale optimization is more difficult than for small-scale studies and we have demonstrated how phantom appearance as well as parameters and convergence of the algorithm can affect recovery studies. For the simplest, and piece-wise constant, structure walnut phantom we found the critical sampling level to be predicted very well by phase-diagram analysis. The situation for the texture walnut phantom was more complex, which motivates further and more extensive large-scale studies, including of the influence of texture on recovery and possibly a different definition of image recovery itself.

%%%%%%%%%%%%%%%%%%%%%%%%%%%%%%%%%%%%%%%%%%%%%%%%%%%%%%%%%%%%%%%%%%%%%
\section{Conclusion} \label{sec:conclusion}

We have presented a systematic framework of phase-diagram analysis from CS for analyzing the undersampling potential of SR in X-ray CT.
In three, quite different, studies we have demonstrated the potential of phase-diagram analysis: We saw that under certain conditions X-ray CT in terms of recoverability performs comparable with an optimal CS sampling strategy of Gaussian sensing matrices; that random sampling in X-ray CT in terms of recoverability does not perform better and in some cases worse than a regular fan-beam sampling setup; and that at least in a simple case the critical sampling level for a large-scale X-ray CT system can be predicted.
An interesting future direction is to address the question: can the observed phase-transition behavior in X-ray CT be theoretically explained, in particular the high degree of similarity with the Gaussian sensing matrix case?

%%%%%%%%%%%%%%%%%%%%%%%%%%%%%%%%%%%%%%%%%%%%%%%%%%%%%%%%%%%%%%%%%%%%%
\section*{Acknowledgment}
The authors are grateful to Martin Lotz for providing code to compute the ALMT phase transitions, to Jared Tanner for providing tabulated values of the DT phase transitions on his website, to Carsten Gundlach for assistance in acquiring the walnut micro-CT dataset, and to Rick Chartrand for inspiring discussions of compressed sensing and tomography. 
The work of JSJ was supported by Advanced Grant 291405 `High-Definition Tomography' from the European Research Council.
This work was supported in part by NIH
Grant Nos. R01-CA158446, R01-CA120540, and R01-EB000225.
The contents of this article are solely the responsibility of
the authors and do not necessarily represent the official
views of the National Institutes of Health.

\bibliographystyle{vancouver}
\bibliography{royal_combined.bib} %cs_bridge

\begin{thebibliography}{10}

\bibitem{SidkyTV:06}
Sidky EY, Kao CM, Pan X.
\newblock {Accurate image reconstruction from few-views and limited-angle data
  in divergent-beam CT}.
\newblock J Xray Sci Technol. 2006;14:119--139.

\bibitem{song2007sparseness}
Song J, Liu QH, Johnson GA, Badea CT.
\newblock {Sparseness prior based iterative image reconstruction for
  retrospectively gated cardiac micro-CT}.
\newblock Med Phys. 2007;34:4476--4483.

\bibitem{sidky2008image}
Sidky EY, Pan X.
\newblock {Image reconstruction in circular cone-beam computed tomography by
  constrained, total-variation minimization}.
\newblock Phys Med Biol. 2008;53:4777--4807.

\bibitem{chen2008prior}
Chen GH, Tang J, Leng S.
\newblock {Prior image constrained compressed sensing (PICCS): A method to
  accurately reconstruct dynamic CT images from highly undersampled projection
  data sets}.
\newblock Med Phys. 2008;35:660--663.

\bibitem{Bian:10}
Bian J, Siewerdsen JH, Han X, Sidky EY, Prince JL, Pelizzari CA, et~al.
\newblock {Evaluation of sparse-view reconstruction from flat-panel-detector
  cone-beam CT}.
\newblock Phys Med Biol. 2010;55:6575--6599.

\bibitem{ritschl2011improved}
Ritschl L, Bergner F, Fleischmann C, Kachelrie\ss{} M.
\newblock {Improved total variation-based CT image reconstruction applied to
  clinical data}.
\newblock Phys Med Biol. 2011;56:1545--1561.

\bibitem{Donoho2006}
Donoho DL.
\newblock {Compressed sensing}.
\newblock IEEE Trans Inf Theory. 2006;52:1289--1306.

\bibitem{CandesTao2005decoding}
Candes EJ, Tao T.
\newblock {Decoding by linear programming}.
\newblock IEEE Trans Inf Theory. 2005;51:4203--4215.

\bibitem{candes2006robust}
Cand\`{e}s EJ, Romberg J, Tao T.
\newblock {Robust uncertainty principles: Exact signal reconstruction from
  highly incomplete frequency information}.
\newblock IEEE Trans Inf Theory. 2006;52:489--509.

\bibitem{delaney1998globally}
Delaney AH, Bresler Y.
\newblock {Globally convergent edge-preserving regularized reconstruction: an
  application to limited-angle tomography}.
\newblock IEEE Trans Image Process. 1998;7:204--221.

\bibitem{li2002accurate}
Li M, Yang H, Kudo H.
\newblock {An accurate iterative reconstruction algorithm for sparse objects:
  application to 3D blood vessel reconstruction from a limited number of
  projections}.
\newblock Phys Med Biol. 2002;47:2599--2609.

\bibitem{candes2008introduction}
Cand\`{e}s EJ, Wakin MB.
\newblock {An introduction to compressive sampling}.
\newblock IEEE Signal Process Mag. 2008;21:21--30.

\bibitem{Tillmann2014}
Tillmann AM, Pfetsch ME.
\newblock {The computational complexity of the restricted isometry property,
  the nullspace property, and related concepts in compressed sensing}.
\newblock IEEE Trans Inf Theory. 2014;60:1248--1259.

\bibitem{CandesTao2006nearoptimal}
Candes EJ, Tao T.
\newblock {Near-optimal signal recovery from random projections: universal
  encoding strategies?}
\newblock IEEE Trans Inf Theory. 2006;52:5406--5425.

\bibitem{FoucartRauhut:2013}
Foucart S, Rauhut H.
\newblock {A Mathematical Introduction to Compressive Sensing}.
\newblock New York, NY: Springer; 2013.

\bibitem{Brady_CStomo:14}
Kaganovsky Y, Li D, Holmgren A, Jeon H, MacCabe KP, Politte DG, et~al.
\newblock {Compressed sampling strategies for tomography}.
\newblock J Opt Soc Amer A. 2014;31:1369--1394.

\bibitem{GilbertIndyk2010}
Gilbert A, Indyk P.
\newblock {Sparse recovery using sparse matrices}.
\newblock Proc IEEE. 2010;98:937--947.

\bibitem{elad2010sparse}
Elad M.
\newblock {Sparse and Redundant Representations: From Theory to Applications in
  Signal and Image Processing}.
\newblock New York, NY: Springer; 2010.

\bibitem{Petra:2009}
Petra S, Schn\"{o}rr C.
\newblock {TomoPIV meets compressed sensing}.
\newblock Pure Math Appl. 2009;20:49--76.

\bibitem{PetraSchnoerr2014}
Petra S, Schn\"{o}rr C.
\newblock {Average case recovery analysis of tomographic compressive sensing}.
\newblock Linear Algebra Appl. 2014;441:168--198.

\bibitem{Joergensen_TMI:2013}
J\o{}rgensen JS, Sidky EY, Pan X.
\newblock {Quantifying admissible undersampling for sparsity-exploiting
  iterative image reconstruction in x-ray CT}.
\newblock IEEE Trans Med Imaging. 2013;32:460--473.

\bibitem{Joergensen_eqconpap_v2_arxiv:2014}
J\o{}rgensen JS, Sidky EY, Hansen PC, Pan X.
\newblock {Empirical average-case relation between undersampling and sparsity
  in x-ray CT}; 2014.
\newblock Available from: \url{http://arxiv.org/abs/1211.5658}.

\bibitem{uniqueness_arxiv_2014}
J\o{}rgensen JS, Kruschel C, Lorenz DA.
\newblock {Testable uniqueness conditions for empirical assessment of
  undersampling levels in total variation-regularized X-ray CT}.
\newblock Inverse Probl Sci Eng, to appear. 2014;Available from:
  \url{http://arxiv.org/abs/1409.0214}.

\bibitem{DonohoTanner_nonneg:2005}
Donoho DL, Tanner J.
\newblock {Sparse nonnegative solution of underdetermined linear equations by
  linear programming.}
\newblock P Natl Acad Sci USA. 2005;102:9446--51.

\bibitem{DonohoTanner:2009}
Donoho D, Tanner J.
\newblock {Observed universality of phase transitions in high-dimensional
  geometry, with implications for modern data analysis and signal processing}.
\newblock Philos Trans R Soc Lond Ser A. 2009;367:4273--4293.

\bibitem{DonohoTanner_radically:2009}
Donoho DL, Tanner J.
\newblock {Counting faces of randomly projected polytopes when the projection
  radically lowers dimension}.
\newblock J Am Math Soc. 2009;22:1--53.

\bibitem{DonohoTanner_finitesize:2010}
Donoho DL, Tanner J.
\newblock {Exponential bounds implying construction of compressed sensing
  matrices, error-correcting codes, and neighborly polytopes by random
  sampling}.
\newblock IEEE Trans Inf Theory. 2010;56:2002--2016.

\bibitem{Vershik1992}
Vershik AM, Sporyshev PV.
\newblock {Asymptotic behavior of the number of faces of random polyhedra and
  the neighborliness problem}.
\newblock Selecta Math Soviet. 1992;11:181--201.

\bibitem{Amelunxen_arxiv:2014}
Amelunxen D, Lotz M, McCoy MB, Tropp JA.
\newblock {Living on the edge: phase transitions in convex programs with random
  data}; 2014.
\newblock Available from: \url{http://arxiv.org/abs/1303.6672v2}.

\bibitem{needell2013cstv}
Needell D, Ward R.
\newblock {Near-optimal compressed sensing guarantees for total variation
  minimization.}
\newblock IEEE Trans Image Process. 2013;22:3941--3949.

\bibitem{NeedellWard:SIR:13}
Needell D, Ward R.
\newblock {Stable image reconstruction using total variation minimization}.
\newblock SIAM J Imaging Sci. 2013;6:1035--1058.

\bibitem{Bayati:2012}
Bayati M, Lelarge M, Montanari A.
\newblock {Universality in polytope phase transitions and message passing
  algorithms}; 2012.
\newblock Available from: \url{http://arxiv.org/abs/1207.7321}.

\bibitem{MOSEK}
{MOSEK ApS}. {MOSEK Optimization Software, version 6.0.0.122}.
\newblock Copenhagen, Denmark; 2011.
\newblock Available from: \url{www.mosek.com}.

\bibitem{Hansen2012}
Hansen PC, Saxild-Hansen M.
\newblock {AIR Tools -- A MATLAB package of algebraic iterative reconstruction
  methods}.
\newblock J Comput Appl Math. 2012;236:2167--2178.

\bibitem{SiltanenTV:2014}
H\"{a}m\"{a}l\"{a}inen K, Harhanen L, Hauptmann A, Kallonen A, Niemi E,
  Siltanen S.
\newblock {Total variation regularization for large-scale X-ray tomography}.
\newblock Int J Tomo Simulat. 2014;25:1--25.

\bibitem{beck2009fast}
Beck A, Teboulle M.
\newblock A fast iterative shrinkage-thresholding algorithm for linear inverse
  problems.
\newblock SIAM J Imag Sci. 2009;2:183--202.

\bibitem{chambolle2011first}
Chambolle A, Pock T.
\newblock A first-order primal-dual algorithm for convex problems with
  applications to imaging.
\newblock J Math Imag Vis. 2011;40:120--145.

\bibitem{jakob:2013}
J{\o}rgensen JS, Sidky EY, Pan X.
\newblock Quantifying admissible undersampling for sparsity-exploiting
  iterative image reconstruction in {X}-ray {CT}.
\newblock IEEE Trans Med Imag. 2013;32:460--473.

\bibitem{sidky:CP:2012}
Sidky EY, J{\o}rgensen JH, Pan X.
\newblock Convex optimization problem prototyping for image reconstruction in
  computed tomography with the {C}hambolle-{P}ock algorithm.
\newblock Phys Med Biol. 2012;57:3065--3091.

\bibitem{sidky2013first}
Sidky EY, J{\o}rgensen JS, Pan X.
\newblock First-order convex feasibility algorithms for {X}-ray {CT}.
\newblock Med Phys. 2013;40(3):031115.

\bibitem{SidkyJTEHM:2014}
Sidky EY, Chartrand R, Boone JM, Pan X.
\newblock Constrained {T$p$V} minimization for enhanced exploitation of
  gradient sparsity: application to {CT} image reconstruction.
\newblock J Trans Eng Health Med. 2014;2:1800418.

\bibitem{de2004distance}
Man BD, Basu S.
\newblock Distance-driven projection and backprojection in three dimensions.
\newblock Phys Med Biol. 2004;49:2463--2475.

\bibitem{zeng2000unmatched}
Zeng GL, Gullberg GT.
\newblock Unmatched projector/backprojector pairs in an iterative
  reconstruction algorithm.
\newblock IEEE Trans Med Imag. 2000;19:548--555.

\end{thebibliography}
 
\end{document}